# INTRINSIC ULTRACONTRACTIVITY OF NONSYMMETRIC DIFFUSIONS WITH MEASURE-VALUED DRIFTS AND POTENTIALS

By Panki Kim[1] and Renming Song[2]

*Seoul National University and University of Illinois*

Recently, in [Preprint (2006)], we extended the concept of intrinsic ultracontractivity to nonsymmetric semigroups. In this paper, we study the intrinsic ultracontractivity of nonsymmetric diffusions with measure-valued drifts and measure-valued potentials in bounded domains. Our process $Y$ is a diffusion process whose generator can be formally written as $L+\mu\cdot\nabla-\nu$ with Dirichlet boundary conditions, where $L$ is a uniformly elliptic second-order differential operator and $\mu=(\mu^1,\dots,\mu^d)$ is such that each component $\mu^i$, $i=1,\dots,d$, is a signed measure belonging to the Kato class $\mathbf{K}_{d,1}$ and $\nu$ is a (nonnegative) measure belonging to the Kato class $\mathbf{K}_{d,2}$. We show that scale-invariant parabolic and elliptic Harnack inequalities are valid for $Y$.

In this paper, we prove the parabolic boundary Harnack principle and the intrinsic ultracontractivity for the killed diffusion $Y^D$ with measure-valued drift and potential when $D$ is one of the following types of bounded domains: twisted Hölder domains of order $\alpha\in(1/3,1]$, uniformly Hölder domains of order $\alpha\in(0,2)$ and domains which can be locally represented as the region above the graph of a function. This extends the results in [*J. Funct. Anal.* **100** (1991) 181–206] and [*Probab. Theory Related Fields* **91** (1992) 405–443]. As a consequence of the intrinsic ultracontractivity, we get that the supremum of the expected conditional lifetimes of $Y^D$ is finite.

**1. Introduction.** In this paper, we study the intrinsic ultracontractivity of a nonsymmetric diffusion process $Y$ with measure-valued drift and

Received May 2006; revised October 2007.

[1]Supported by the Korea Research Foundation Grant funded by the Korean Government (MOEHRD, Basic Research Promotion Fund) (KRF-2007-331-C00037).

[2]Supported in part by a joint US-Croatia Grant INT 0302167.

measure-valued potential in bounded domains $D \subset \mathbf{R}^d$ for $d \geq 3$. The generator of $Y$ can be formally written as $L + \mu \cdot \nabla - \nu$ with Dirichlet boundary conditions, where $L$ is a uniformly elliptic second order differential operator and $\mu = (\mu^1, \dots, \mu^d)$ is such that each component $\mu^i$, $i = 1, \dots, d$, is a signed measure belonging to the Kato class $\mathbf{K}_{d,1}$ and $\nu$ is a (nonnegative) measure belonging to the Kato class $\mathbf{K}_{d,2}$ (see below for the definitions of $\mathbf{K}_{d,1}$ and $\mathbf{K}_{d,2}$). The existence and uniqueness of this process $Y$ were proven in Bass and Chen [3]. In [15, 16, 17, 19], we have studied properties of diffusions with measure-valued drifts in bounded domains. Using results in [15, 16, 17, 19], we will prove that, with respect to a certain reference measure, $Y$ has a dual process which is a continuous Hunt process satisfying the strong Feller property.

The notion of intrinsic ultracontractivity, introduced in [11] for symmetric semigroups, is a very important concept and has been studied extensively. In [18], the concept of intrinsic ultracontractivity was extended to nonsymmetric semigroups and it was proven there that the semigroup of the killed diffusion process in a bounded Lipschitz domain is intrinsically ultracontractive if the coefficients of the generator of the diffusion process are smooth.

In this paper, using the duality of our processes, we prove that the semigroups of the killed diffusion $Y^D$ and its dual are intrinsically ultracontractive if $D$ is one of the following types of bounded domains:

(a) a twisted Hölder domain of order $\alpha \in (1/3, 1]$;
(b) a uniformly Hölder domain of order $\alpha \in (0, 2)$;
(c) a domain which can be locally represented as the region above the graph of a function.

In fact, we first prove parabolic boundary Harnack principles for $Y^D$ and its dual process (see Theorem 5.6 and Corollary 5.7). We then show that the parabolic boundary Harnack principles imply that the semigroups of $Y^D$ and its dual are intrinsically ultracontractive. The fact that the parabolic boundary Harnack principle implies the intrinsic ultracontractivity in the symmetric diffusion case was used and discussed in [2] and [12]. As a consequence of the intrinsic ultracontractivity, we have that the supremum of the expected conditional lifetimes of $Y^D$ is finite if $D$ is one of the domains above.

Many results in this paper are stated for both the diffusion process $Y$ and its dual. In these cases, the proofs for the dual process are usually harder. Once the proofs for the dual process are completed, it is very easy to see that the results for the diffusion process $Y$ can be proven through similar and simpler arguments. For this reason, we only present the proof for the dual process.

The contents of this paper are organized as follows. In Section 2, we present some preliminary properties of $Y$ and the existence of the dual process of $Y$. Section 3 contains the proof of parabolic Harnack inequalities for $Y$ and its dual process. In Section 4, we discuss some properties of $Y$ and its dual in twisted Hölder domains, uniformly Hölder domains and domains which can be locally represented as the region above the graph of a function. In the last section, we prove the parabolic boundary Harnack principles and show that the parabolic boundary Harnack principles imply the intrinsic ultracontractivity of the nonsymmetric semigroups. Finally, we obtain that the supremum of the expected conditional lifetime is finite.

In this paper, we always assume that $d \geq 3$. Throughout, we use the notation $a \wedge b := \min\{a,b\}$ and $a \vee b := \max\{a,b\}$. The distance between $x$ and $\partial D$, the boundary of $D$, is denote by $\rho_D(x)$. We use the convention $f(\partial) = 0$. We will also use the following convention: the values of the constants $r_1, t_0, t_1$ will remain the same throughout, while the values of the constants $c_1, c_2, \ldots$ might change from one appearance to another. The labeling of the constants $c_1, c_2, \ldots$ starts anew in the statement of each result.

In this paper, we use ":=" to denote a definition, this being read as "is defined to be."

## 2. Dual processes for diffusion processes with measure-valued drifts and potentials.

First, we recall the definition of the Kato class $\mathbf{K}_{d,j}$ for $j = 1, 2$. For any function $f$ on $\mathbf{R}^d$ and $r > 0$, we define

$$M_f^j(r) = \sup_{x \in \mathbf{R}^d} \int_{|x-y| \leq r} \frac{|f|(y)\, dy}{|x-y|^{d-j}}, \qquad j = 1, 2.$$

For any signed measure $\nu$ on $\mathbf{R}^d$, we use $\nu^+$ and $\nu^-$ to denote its positive and negative parts, and $|\nu| := \nu^+ + \nu^-$. For any signed measure $\nu$ on $\mathbf{R}^d$ and any $r > 0$, we define

$$M_\nu^j(r) = \sup_{x \in \mathbf{R}^d} \int_{|x-y| \leq r} \frac{|\nu|(dy)}{|x-y|^{d-j}}, \qquad j = 1, 2.$$

DEFINITION 2.1. Let $j = 1, 2$. We say that a function $f$ on $\mathbf{R}^d$ *belongs to the Kato class* $\mathbf{K}_{d,j}$ if $\lim_{r \downarrow 0} M_f^j(r) = 0$. We say that a signed Radon measure $\nu$ on $\mathbf{R}^d$ *belongs to the Kato class* $\mathbf{K}_{d,j}$ if $\lim_{r \downarrow 0} M_\nu^j(r) = 0$.

Throughout this paper, we assume that $\mu = (\mu^1, \ldots, \mu^d)$ and $\nu$ are fixed with each $\mu^i$ being a signed measure on $\mathbf{R}^d$ belonging to $\mathbf{K}_{d,1}$ and $\nu$ being a (non-negative) measure on $\mathbf{R}^d$ belonging to $\mathbf{K}_{d,2}$.

We also assume that the operator $L$ is either $L_1$ or $L_2$, where

$$L_1 := \frac{1}{2} \sum_{i,j=1}^d \partial_i (a_{ij} \partial_j) \quad \text{and} \quad L_2 := \frac{1}{2} \sum_{i,j=1}^d a_{ij} \partial_i \partial_j$$

with $\mathbf{A}(x) := (a_{ij}(x))$ being $C^1$ and uniformly elliptic. Since $\mathbf{A}(x) = (a_{ij}(x))$ is $C^1$, without loss of generality, one can assume that the matrix $A(x)$ is symmetric (see, e.g., Section 6 of [16]).

We will use $X$ to denote the diffusion process in $\mathbf{R}^d$ whose generator can be formally written as $L + \mu \cdot \nabla$. When each $\mu^i$ is given by $U^i(x)\, dx$ for some function $U^i$, $X$ is a diffusion in $\mathbf{R}^d$ with generator $L + U \cdot \nabla$ and it is a solution to the stochastic differential equation $dX_t = dX^0_t + U(X_t) \cdot dt$, where $X^0$ is a diffusion in $\mathbf{R}^d$ with generator $L$. For a precise definition of a (nonsymmetric) diffusion $X$ with drift $\mu$ in $\mathbf{K}_{d,1}$, we refer to Section 6 in [16] and Section 1 in [19]. The existence and uniqueness of $X$ were established in [3] (see Remark 6.1 in [3]).

For any open set $U$, we use $\tau^X_U$ to denote the first exit time of $U$ for $X$, that is, $\tau^X_U = \inf\{t > 0 : X_t \notin U\}$. We define $X^U_t(\omega) = X_t(\omega)$ if $t < \tau^X_U(\omega)$ and $X^U_t(\omega) = \partial$ if $t \geq \tau^X_U(\omega)$, where $\partial$ is a cemetery state. The process $X^U$ is called a *killed diffusion with drift $\mu$ in $U$*. $X^U$ is a Hunt process with the strong Feller property, that is, for every $f \in L^\infty(U)$, $\mathbf{E}_x[f(X^U_t)]$ is in $C(U)$, the space of continuous functions in $U$ (Proposition 2.1 [19]). Moreover, $X^U$ has a jointly continuous density $q^U(t, x, y)$ with respect to the Lebesgue measure (Theorem 2.4 in [16]).

From Section 3 in [17] and Proposition 7.1 in [19], we know that for every bounded domain $U$, there exists a positive continuous additive functional $A^U$ of $X^U$ with Revuz measure $\nu|_U$, that is, for any $x \in U$, $t > 0$ and bounded nonnegative function $f$ on $U$,

$$\mathbf{E}_x \int_0^t f(X^U_s)\, dA^U_s = \int_0^t \int_U q^U(s, x, y) f(y) \nu(dy)\, ds.$$

Throughout this paper, we assume that $V$ is a bounded smooth domain in $\mathbf{R}^d$ and consider the transient diffusion process $Y$ such that

$$\mathbf{E}_x[f(Y_t)] = \mathbf{E}_x[\exp(-A^V_t) f(X^V_t)].$$

(See III.3 of [4] for the construction of such a killed process.) We will use $\zeta$ to denote the lifetime of $Y$. Note that the process $Y$ might have killing inside $V$, that is, $\mathbf{P}_x(Y_{\zeta-} \in V)$ might be positive.

A simple example of $Y$ is a diffusion whose infinitesimal generator is a second-order differential operator $L - b \cdot \nabla - c$, where $(b_1, \dots, b_d)$ and $c \geq 0$ belong to the Kato classes $\mathbf{K}_{d,1}$ and $\mathbf{K}_{d,2}$, respectively. If $(b_1, \dots, b_d)$ is differentiable and $L = L_1$, then the formal adjoint of the above operator is $L_1 + b \cdot \nabla - (c - \nabla b)$. If one further assumes that $c - \nabla b \geq 0$, then there is a diffusion process with generator $L_1 + b \cdot \nabla - (c - \nabla b)$. We cannot, and do not make such assumptions in this paper. Instead, we will introduce a new reference measure and consider a dual process with respect to this reference measure.

Recall that for any domain $D\subset \mathbf{R}^d$, $\rho_D(x)$ is the distance between $x$ and $\partial D$. It is shown in [17] that the process $Y$ has a jointly continuous and strictly positive transition density function $r(t,x,y)$ with respect to the Lebesgue measure and for each $T>0$, there exist positive constants $c_j$, $1\le j\le 4$, depending on $V$ such that for $t\le T$,

$$\begin{aligned}
&c_1 t^{-d/2}\left(1\wedge \frac{\rho_V(x)}{\sqrt{t}}\right)\left(1\wedge \frac{\rho_V(y)}{\sqrt{t}}\right)e^{-c_2|x-y|^2/(2t)}\\
&\quad\le r(t,x,y)\\
&\quad\le c_3 t^{-d/2}\left(1\wedge \frac{\rho_V(x)}{\sqrt{t}}\right)\left(1\wedge \frac{\rho_V(y)}{\sqrt{t}}\right)e^{-c_4|x-y|^2/(2t)}.
\end{aligned}\tag{2.1}$$

Moreover, for every smooth subset $U$ of $V$, the killed process $Y^U$ has a jointly continuous and strictly positive transition density function $r^U(t,x,y)$ with respect to the Lebesgue measure and for each $T>0$, there exist positive constants $c_j$, $5\le j\le 8$, depending on $U$ such that for $t\le T$,

$$\begin{aligned}
&c_5 t^{-d/2}\left(1\wedge \frac{\rho_U(x)}{\sqrt{t}}\right)\left(1\wedge \frac{\rho_U(y)}{\sqrt{t}}\right)e^{-c_6|x-y|^2/(2t)}\\
&\quad\le r^U(t,x,y)\\
&\quad\le c_7 t^{-d/2}\left(1\wedge \frac{\rho_U(x)}{\sqrt{t}}\right)\left(1\wedge \frac{\rho_V(y)}{\sqrt{t}}\right)e^{-c_8|x-y|^2/(2t)}.
\end{aligned}\tag{2.2}$$

(See Theorem 4.4(1) in [17].)

Let $C_0(V)$ be the class of bounded continuous functions on $V$ vanishing continuously near the boundary of $V$. We will use $\|\cdot\|_\infty$ to denote the $L^\infty$-norm in $C_0(V)$. Using the joint continuity of $r(t,x,y)$ and $r^U(t,x,y)$, and the estimates above, it is easy to show the following result, so we omit the proof.

PROPOSITION 2.2. *$Y$ is a doubly Feller process (a Feller process satisfying the strong Feller property), that is, for every $g\in C_0(V)$, $\mathbf{E}_x[g(Y_t)]=\mathbf{E}_x[g(Y_t);t<\zeta]$ is in $C_0(V)$ and $\|\mathbf{E}_x[g(Y_t)]-g(x)\|_\infty\to 0$ as $t\to 0$, and for every $f\in L^\infty(V)$, $\mathbf{E}_x[f(Y_t)]$ is bounded and continuous in $V$.*

In particular, the above proposition implies that for any domain $U\subset V$, $Y^U$ is Hunt process with the strong Feller property (see, e.g., [7]).

We will use $G(x,y)$ to denote the Green function of $Y$. For any domain $U\subset V$, we will use $G_U(x,y)$ to denote the Green function of $Y^U$. Thus,

$$\mathbf{E}_x\int_0^\infty f(Y_t)\,dt=\mathbf{E}_x\int_0^\zeta f(Y_t)\,dt=\int_V G(x,y)f(y)\,dy$$

and

$$\mathbf{E}_x \int_0^\infty f(Y_t^U)\, dt = \mathbf{E}_x \int_0^{\tau_U} f(Y_t^U)\, dt = \int_U G_U(x,y) f(y)\, dy,$$

where $\tau_U$ is the first exit time of $U$ for $Y$, that is, $\tau_U = \inf\{t > 0 : Y_t \notin U\}$. We will use $G_U^X(x,y)$ to denote the Green function of $X^U$ and $G_V^0(x,y)$ the Green function of the killed Brownian motion in $V$. Since $Y$ is transient, combining Theorem 6.2 in [15] and the result in Section 3 of [17], we have that there exists a constant $c = c(V)$ such that

$$c^{-1} G_V^0(x,y) \le G(x,y) \le c G_V^0(x,y), \qquad V \times V \setminus \{x = y\}. \tag{2.3}$$

Thus, for every $U \subset V$,

$$G_U(x,y) \le G(x,y) \le \frac{c}{|x-y|^{d-2}} \qquad \text{for every } x, y \in D, \tag{2.4}$$

for some constant $c > 0$.

Let

$$H(x) := \int_V G(y,x)\, dy \quad \text{and} \quad \xi(dx) := H(x)\, dx.$$

It is then easy to check (see the proof of Proposition 2.2 in [19]) that $\xi$ is an excessive measure with respect to $Y$, that is, for every Borel function $f \ge 0$,

$$\int_V f(x) \xi(dx) \ge \int_V \mathbf{E}_x[f(Y_t)] \xi(dx).$$

We define a new transition density function with respect to the reference measure $\xi$ by

$$\overline{r}(t,x,y) := \frac{r(t,x,y)}{H(y)}.$$

Then

$$\overline{G}(x,y) := \int_0^\infty \overline{r}(t,x,y)\, dt = \frac{G(x,y)}{H(y)}$$

is the Green function of $Y$ with respect to the reference measure $\xi(dy)$.

Before we discuss properties of $Y$ any further, we recall some definitions. Recall that $\tau_A = \inf\{t > 0 : Y_t \notin A\}$.

Definition 2.3. Suppose $U$ is an open subset of $V$. A nonnegative Borel function $u$ defined on $U$ is said to be:

(1) *harmonic* with respect to $Y$ in $U$ if

$$u(x) = \mathbf{E}_x[u(Y_{\tau_B})] = \mathbf{E}_x[u(Y_{\tau_B}); \tau_B < \zeta], \qquad x \in B, \tag{2.5}$$

for every bounded open set $B$ with $\overline{B} \subset U$;

(2) *superharmonic* with respect to $Y^U$ if

$$u(x) \geq \mathbf{E}_x[u(Y^U_{\tau_B})], \qquad x \in B,$$

for every bounded open set $B$ with $\overline{B} \subset U$;

(3) *excessive* for $Y^U$ if

$$u(x) \geq \mathbf{E}_x[u(Y^U_t)] = \mathbf{E}_x[u(Y^U_t); t < \zeta], \qquad t > 0, x \in U$$

and

$$u(x) = \lim_{t \downarrow 0} \mathbf{E}_x[u(Y^U_t)], \qquad x \in U;$$

(4) a *potential* for $Y^U$ if it is excessive for $Y^U$ and for every sequence $\{U_n\}_{n \geq 1}$ of open sets with $\overline{U_n} \subset U_{n+1}$ and $\bigcup_n U_n = U$,

$$\lim_{n \to \infty} \mathbf{E}_x[u(Y^U_{\tau_{U_n}})] = 0, \qquad \xi\text{-a.e. } x \in U.$$

A Borel function $u$ defined on $\overline{U}$ is said to be *regular harmonic* with respect to $Y$ in $U$ if $u$ is harmonic with respect to $Y$ in $U$ and (2.5) is true for $B = U$.

Since $Y^U$ is a Hunt processes with the strong Feller property, it is easy to check that $u$ is excessive for $Y^U$ if and only if $f$ is lower-semicontinuous in $U$ and superharmonic with respect to $Y^U$. (See Theorem 4.5.3 in [10] for the Brownian motion case; the proof there can adapted easily to the present case.)

Using (2.1)–(2.2) and the joint continuity of $r(t,x,y)$ and $r^U(t,x,y)$, one can easily check that $G_U(x,y)$ is strictly positive and jointly continuous on $(U \times U) \setminus \{(x,y) : x = y\}$. $G_U(x,y)$ is infinite if and only if $x = y$ (see the proof of Theorem 2.6 in [16]). Thus, by (2.3), we see that $H$ is a strictly positive, bounded continuous function on $V$. Moreover, using the estimates for $G^0_V(x,y)$, one can check that there exists a constant $c = c(V)$ such that

$$c^{-1} \rho_V(x) \leq H(x) \leq c \rho_V(x). \tag{2.6}$$

(See Lemma 6.4 in [19] and its proof.) Now, using the above properties and (2.4), we see that $Y$ is a transient diffusion with its Green function $\overline{G}(x,y)$ with respect to $\xi$ satisfying the conditions in [9] and [23] (see (A1)–(A4) in [19]). In fact, one can follow the arguments in [19] and check that all the results in Sections 2–3 of [19] are true for $Y$. In particular, using the same arguments in the proofs of Theorems 2.4–2.5 in [19], it is easy to check that the conditions (i)–(vii) and (70)–(71) in [20] (also, see the Remark on page 391 in [21]) are satisfied. Thus, with respect to the reference measure $\xi$, $Y$ has a nice dual process. For more detailed arguments, we refer readers to [19].

THEOREM 2.4. *There exists a continuous transient Hunt process $\widehat{Y}$ in $V$ such that $\widehat{Y}$ is a strong dual of $Y$ with respect to the measure $\xi$, that is, the density of the semigroup $\{\widehat{P}_t\}_{t\geq 0}$ of $\widehat{Y}$ is $\widehat{r}(t,x,y) := \overline{r}(t,y,x)$ and thus*

$$\int_V f(x)P_t g(x)\xi(dx) = \int_V g(x)\widehat{P}_t f(x)\xi(dx) \qquad \text{for all } f,g \in L^2(V,\xi).$$

We will use $\widehat{\zeta}$ to denote the lifetime of $\widehat{Y}$. Note that $\widehat{Y}$ might also have killing inside $V$, that is, $\mathbf{P}_x(\widehat{Y}_{\widehat{\zeta}-} \in V)$ might be positive.

By Theorem 2 and Remark 2 following it in [25], for any domain $U \subset V$, $Y^U$ and $\widehat{Y}^U$ are duals of each other with respect to $\xi$. For any domain $U \subset V$, we define

$$\widehat{r}^U(t,x,y) := \frac{r^U(t,y,x)H(y)}{H(x)}.$$

Since $H$ is strictly positive and continuous, by the joint continuity of $r^U(t,x,y)$ (see Section 4 of [17] and the references therein), $\widehat{r}^U(t,x,y)$ is jointly continuous on $U \times U$. Thus, $\widehat{r}^U(t,x,y)$ is the transition density of $\widehat{Y}^U$ with respect to the Lebesgue measure and

$$\widehat{G}_U(x,y) := \frac{G_U(y,x)H(y)}{H(x)} \tag{2.7}$$

is the Green function for $\widehat{Y}^U$ with respect to the Lebesgue measure so that for every nonnegative Borel function $f$,

$$\mathbf{E}_x\left[\int_0^{\widehat{\tau}_U} f(\widehat{Y}_t)\,dt\right] = \int_U \widehat{G}_U(x,y)f(y)\,dy,$$

where $\widehat{\tau}_U := \inf\{t > 0 : \widehat{Y}_t \notin U\}$.

We will use $\{\widehat{G}^U_\lambda, \lambda \geq 0\}$ to denote the resolvent of $\widehat{Y}^U$ with respect to $\xi$. Following the argument in Proposition 3.4 in [19], one can check that $\widehat{Y}^U$ has the strong Feller property. We include the proof here for the reader's convenience.

PROPOSITION 2.5. *For any $U \subset V$, $\widehat{Y}^U$ has the strong Feller property in the resolvent sense; that is, for every bounded Borel function $f$ on $U$ and $\lambda \geq 0$, $\widehat{G}^U_\lambda f(x)$ is a bounded continuous function on $U$.*

PROOF. By the resolvent equation $\widehat{G}^U_0 = \widehat{G}^U_\lambda + \lambda \widehat{G}^U_0 \widehat{G}^U_\lambda$, it is enough to show the strong Feller property for $\widehat{G}^U_0$. Fix a bounded Borel function $f$ on $U$ and a sequence $\{y_n\}_{n\geq 1}$ converging to $y$ in $U$. Let $M := \|fH\|_{L_\infty(U)} < \infty$. We assume $\{y_n\}_{n\geq 1} \subset K$ for a compact subset $K$ of $U$. Let $A := \inf_{y\in K} H(y)$. By (2.6), we know that $A$ is strictly positive. Note that there exists a constant $c_1$ such that for every $\delta > 0$,

$$\left(\int_{B(y,\delta)} \frac{dx}{|x-y|^{d-2}} + \int_{B(y_n,2\delta)} \frac{dx}{|x-y_n|^{d-2}}\right) \le c_1\delta^2.$$

Thus, by (2.4), there exists a constant $c_2$ such that for every $\delta > 0$ and $y_n$ with $y_n \in B(y, \frac{\delta}{2}) \subset B(y, 2\delta) \in K$,

$$\begin{aligned} &\int_{B(y,\delta)} \frac{G_U(x,y)H(x)f(x)}{H(y)}\,dx + \int_{B(y,\delta)} \frac{G_U(x,y_n)H(x)f(x)}{H(y_n)}\,dx \\ &\quad \le \frac{M}{A}\left(\int_{B(y,\delta)} G_U(x,y)\,dx + \int_{B(y_n,2\delta)} G_U(x,y_n)\,dx\right) \\ &\quad \le \frac{c_2 M}{A}\left(\int_{B(y,\delta)} \frac{dx}{|x-y|^{d-2}} + \int_{B(y_n,2\delta)} \frac{dx}{|x-y_n|^{d-2}}\right) \le \frac{1}{A} c_1 c_2 M \delta^2. \end{aligned}$$

Given $\varepsilon$, choose $\delta$ small enough such that $\frac{1}{A} c_1 c_2 M\delta^2 < \frac{\varepsilon}{2}$. Then

$$\begin{aligned} &|\widehat{G}_0^U f(y) - \widehat{G}_0^U f(y_n)| \\ &\qquad \le M \int_{U\setminus B(y,\delta)} \left|\frac{G_U(x,y)}{H(y)} - \frac{G_U(x,y_n)}{H(y_n)}\right| dx + \frac{\varepsilon}{2}. \end{aligned}$$

Note that $G_U(x,y_n)/H(y_n)$ converges to $G_U(x,y)/H(y)$ for every $x \ne y$ and that $\{G_U(x,y_n)/H(y_n)\}$ are uniformly bounded on $x \in U \setminus B(y,\delta)$ and $y_n \in B(y, \frac{\delta}{2})$. So, the first term on the right-hand side of the inequality above goes to zero as $n \to \infty$, by the bounded convergence theorem. □

Applying the results in [23] and [24], we have the following.

Proposition 2.6. *Suppose $D \subset V$. Any function which is harmonic for $Y$ (resp. $\widehat{Y}$) in $D$ is continuous. For each $y$, $x \to G_D(x,y)$ is excessive for $Y^D$ and harmonic for $Y$ in $D \setminus \{y\}$, and $x \to \widehat{G}_D(x,y)$ is excessive for $\widehat{Y}^D$ and harmonic for $\widehat{Y}$ in $D \setminus \{y\}$. Moreover, for every open subset $U$ of $D$, we have*

$$\begin{aligned} &\mathbf{E}_x[G_D(Y^D_{T_U}, y)] = G_D(x,y) \quad \textit{and} \\ &\mathbf{E}_x[\widehat{G}_D(\widehat{Y}^D_{\widehat{T}_U}, y)] = \widehat{G}_D(x,y), (x,y) \in D \times U, \end{aligned} \tag{2.8}$$

*where $T_U := \inf\{t > 0 : Y^D_t \in U\}$ and $\widehat{T}_U := \inf\{t > 0 : \widehat{Y}^D_t \in U\}$. In particular, for every $y \in D$ and $\varepsilon > 0$, $G_D(\cdot, y)$ is regular harmonic with respect to $Y^D$ in $D \setminus B(y,\varepsilon)$ and $\widehat{G}_D(\cdot,y)$ is regular harmonic with respect to $\widehat{Y}^D$ in $D \setminus B(y,\varepsilon)$.*

By Theorem 3.7 in [16], there exist constants $r_1 = r_1(d, \mu) > 0$ and $c = c(d, \mu) > 1$ depending on $\mu$ only via the rate at which $\max_{1\le i\le d} M^1_{\mu^i}(r)$ goes to zero such that for $r \le r_1$, $w \in \mathbf{R}^d$, $x, y \in B(w, r)$,

$$c^{-1} G^0_{B(w,r)}(x, y) \le G^X_{B(w,r)}(x, y) \le c G^0_{B(w,r)}(x, y). \tag{2.9}$$

Thus, there exists a positive constant $c$ independent of $r \le r_1$ such that for every $x, y, z \in B(w, r)$ and $w \in \mathbf{R}^d$,

$$\frac{G^X_{B(w,r)}(x, y) G^X_{B(w,r)}(y, z)}{G^X_{B(w,r)}(x, z)} \le c(|x - y|^{2-d} + |y - z|^{2-d}). \tag{2.10}$$

For any $z \in B(w, r)$, let $(\mathbf{P}^z_x, X^{B(w,r)}_t)$ be the $G^X_{B(w,r)}(\cdot, z)$-transform of $(\mathbf{P}_x, X^{B(w,r)}_t)$, that is, for any nonnegative Borel function $f$,

$$\mathbf{E}^z_x[f(X^{B(w,r)}_t)] = \mathbf{E}_x\Big[\frac{G^X_{B(w,r)}(X^{B(w,r)}_t, z)}{G^X_{B(w,r)}(x, z)} f(X^{B(w,r)}_t)\Big].$$

Recall that $A^V$ is the positive continuous additive functionals of $X^V$ with Revuz measures $\nu|_V$. Equation (2.10) implies that there exists a positive constant $c_1 < \infty$ such that for every $r \in (0, r_1]$, $w \in \mathbf{R}^d$ and $x, z \in B(w, r)$,

$$\mathbf{E}^z_x[A^V_{\tau^X_{B(w,r)}}] \le \int_{B(w,r)} \frac{G^X_{B(w,r)}(x, y) G^X_{B(w,r)}(y, z)}{G^X_{B(w,r)}(x, z)} \nu(dy) < c_1. \tag{2.11}$$

Hence, by Jensen's inequality, for $x, z \in B(w, r)$, we have

$$\mathbf{E}^z_x[\exp(-A^V_{\tau^X_{B(w,r)}})] \ge \exp(-\mathbf{E}^z_x[A^V_{\tau^X_{B(w,r)}}]) \ge e^{-c_1} > 0.$$

Combining the identity

$$G_{B(w,r)}(x, z) = G^X_{B(w,r)}(x, z) \mathbf{E}^z_x[\exp(-A^V_{\tau^X_{B(w,r)}})], \qquad x, z \in B(w, r),$$

(Lemma 3.5 (1) of [5]) with (2.9), we arrive at the following result.

Proposition 2.7. *There exist positive constants $c$ and $r_1 := r_1(d, \mu, \nu)$ such that for all $r \in (0, r_1]$ and $B(w, r) \in V$, we have*

$$c^{-1} G^0_{B(w,r)}(x, y) \le G_{B(w,r)}(x, y) \le c G^0_{B(w,r)}(x, y), \qquad x, y \in B(w, r).$$

In the remainder of this paper, we will always assume $D$ is a bounded domain with $\overline{D} \subset V$. Let $\gamma_1 := \frac{1}{2}\mathrm{dist}(\partial V, \overline{D})$ and $\check{V} := \{z \in V; \rho_V(z) > \gamma_1\}$.

We fix $D$, $\check{V}$ and $\gamma_1$ throughout this paper. For any subdomain $U \subset V$ and any subset $A$ of $U$, we define

$$\begin{aligned}\mathrm{Cap}^U(A) := \sup\biggl\{\eta(A) : \eta \text{ is a measure supported on } A \\ \text{with } \int_U G_U^0(x,y)\eta(dy) \le 1\biggr\}.\end{aligned} \tag{2.12}$$

The next lemma is a nonsymmetric version of Lemma 2.1 in [2] for small balls. For any set $A$, we define $A_r^z := z + rA = \{w \in \mathbf{R}^d : w = z + ra, a \in A\}$, $A_r := A_r^0$ and $A^z := A_1^z$.

Lemma 2.8. *There exists $c = c(V, d, \mu, \nu) > 0$ such that for any compact subset $K$ of $B(0,1)$, $r \in (0, r_1]$, $B(z,r) \subset \check{V}$ and compact set $A \subset K_r$, we have for any $x \in B(z,r)$,*

$$\begin{aligned}&c^{-1} r^{2-d}\Bigl(\inf_{y\in K} G^0_{B(0,1)}((x-z)/r, y)\Bigr)\mathrm{Cap}^{B(0,r)}(A) \\ &\quad \le \mathbf{P}_x(T_{A^z} < \tau_{B(z,r)}) \\ &\quad \le c r^{2-d}\Bigl(\sup_{y\in K} G^0_{B(0,1)}((x-z)/r, y)\Bigr)\mathrm{Cap}^{B(0,r)}(A)\end{aligned}$$

*and*

$$\begin{aligned}&c^{-1} r^{2-d}\Bigl(\inf_{y\in K} G^0_{B(0,1)}((x-z)/r, y)\Bigr)\mathrm{Cap}^{B(0,r)}(A) \\ &\quad \le \mathbf{P}_x(\widehat{T}_{A^z} < \widehat{\tau}_{B(z,r)}) \\ &\quad \le c r^{2-d}\Bigl(\sup_{y\in K} G^0_{B(0,1)}((x-z)/r, y)\Bigr)\mathrm{Cap}^{B(0,r)}(A).\end{aligned}$$

Proof. For $B(z,r) \subset \check{V}$ and $U \subset B(z,r)$, define

$$\begin{aligned}\mathrm{Cap}_{\widehat{Y}}^{B(z,r)}(U) := \sup\biggl\{\eta(U) : \eta \text{ is a measure supported on } U \\ \text{with } \int_{B(z,r)} \widehat{G}_{B(z,r)}(x,y)\eta(dy) \le 1\biggr\}.\end{aligned} \tag{2.13}$$

From (2.6), (2.7) and Proposition 2.7, we see that there is a constant $c > 0$ such that for every $r < r_1$ and $B(z,r) \subset \check{V}$, we have

$$c^{-1}\mathrm{Cap}_{\widehat{Y}}^{B(z,r)}(U) \le \mathrm{Cap}^{B(z,r)}(U) \le c\,\mathrm{Cap}_{\widehat{Y}}^{B(z,r)}(U), \qquad U \subset B(z,r). \tag{2.14}$$

Note that $Y^{B(z,r)}$ and $\widehat{Y}^{B(z,r)}$ are Hunt processes with the strong Feller property and that they are in the strong duality with respect to $\xi$ (Propositions 2.4 and 2.5). Since $A^z$ is a compact subset of $B(z,r)$, there exist capacitary measures $\mu_{A^z}$ for $A^z$ with respect to $Y^{B(z,r)}$ and $\widehat{\mu}_{A^z}$ for $A^z$ with respect to $\widehat{Y}^{B(z,r)}$ such that $\mathrm{Cap}^{B(z,r)}_{\widehat{Y}}(A^z) = \mu_{A^z}(A^z) = \widehat{\mu}_{A^z}(A^z)$. (See, e.g., VI.4 of [4] and Sections 5.1–5.2 of [10] for details.)

Using Proposition 2.7 and (2.6), we have for every $x \in B(z,r)$,

$$\begin{aligned} \int_{A^z} \widehat{G}_{B(z,r)}(x,y)\widehat{\mu}_{A^z}(dy) &= \int_{A^z} \frac{G_{B(z,r)}(y,x)H(y)}{H(x)}\widehat{\mu}_{A^z}(dy) \\ &\ge c_1^{-1}\int_{A^z} G^0_{B(z,r)}(x,y)\widehat{\mu}_{A^z}(dy) \\ &\ge c_1^{-1}\Big(\inf_{y\in K^z_r} G^0_{B(z,r)}(x,y)\Big)\widehat{\mu}_{A^z}(A^z) \\ &= c_1^{-1}\Big(\inf_{y\in K^z_r} G^0_{B(z,r)}(x,y)\Big)\,\mathrm{Cap}^{B(z,r)}_{\widehat{Y}}(A^z) \end{aligned} \tag{2.15}$$

for some constant $c_1 > 0$. Applying (2.14) to the above equation and using the scaling property of Brownian motion, we get that for every $x \in B(z,r)$,

$$\begin{aligned} &\Big(\inf_{y\in K^z_r} G^0_{B(z,r)}(x,y)\Big)\,\mathrm{Cap}^{B(z,r)}_{\widehat{Y}}(A^z) \\ &\qquad \ge c^{-1}r^{2-d}\Big(\inf_{y\in K} G^0_{B(0,1)}((x-z)/r,y)\Big)\,\mathrm{Cap}^{B(0,r)}(A). \end{aligned} \tag{2.16}$$

On the other hand, by (2.8), we have for every $x \in B(z,r)$,

$$\begin{aligned} &\int_{A^z} \widehat{G}_{B(z,r)}(x,y)\widehat{\mu}_{A^z}(dy) \\ &\qquad = \int_{A^z} \mathbf{E}_x[\widehat{G}_{B(z,r)}(\widehat{Y}^{B(z,r)}_{\widehat{T}_{A^z}},y)]\ \widehat{\mu}_{A^z}(dy) \\ &\qquad \le \Big(\sup_{w\in A^z}\int_{A^z} \widehat{G}_{B(z,r)}(w,y)\mu_{A^z}(dy)\Big)\mathbf{P}_x(\widehat{T}_{A^z} < \widehat{\tau}_{B(z,r)}) \\ &\qquad \le c_2\mathbf{P}_x(\widehat{T}_{A^z} < \widehat{\tau}_{B(z,r)}) \end{aligned} \tag{2.17}$$

for some constant $c_2 > 0$. In the last inequality above, we have used (2.6) and (2.13).

Combining (2.15)–(2.17), we have for every $x \in B(z,r)$,

$$\mathbf{P}_x(\widehat{T}_{A^z} < \widehat{\tau}_{B(z,r)}) \ge c_3 r^{2-d}\Big(\inf_{y\in K} G^0_{B(0,1)}((x-z)/r,y)\Big)\mathrm{Cap}^{B(0,r)}(A)$$

for some constant $c_3 > 0$. Thus, we have shown the first inequality in (2).

By Corollary 1 to Theorem 2 in [9], the function $x \mapsto \mathbf{P}_x(\widehat{T}_{A^z} < \widehat{\tau}_{B(z,r)})$ is a potential for $\widehat{Y}^{B(z,r)}$, thus there exists a Radon measure $\widehat{\nu}_1$ on $A^z$ such

that

$$\mathbf{P}_x(\widehat{T}_{A^z} < \widehat{\tau}_{B(z,r)}) = \int_{A^z} \widehat{G}_{B(z,r)}(x,y)\widehat{\nu}_1(dy), \qquad x \in B(z,r).$$

Hence, by (2.6) and (2.13), we have

$$\mathbf{P}_x(\widehat{T}_{A^z} < \widehat{\tau}_{B(z,r)}) \le c_4\Big(\sup_{y\in K_r^z} G_{B(z,r)}(y,x)\Big)\,\mathrm{Cap}_{\widehat{Y}}^{B(z,r)}(A^z), \qquad x \in B(z,r)$$

for some constant $c_4 > 0$. Now, applying Proposition 2.7 and (2.14) to the right-hand side above and using the scaling property of Brownian motion, we get the desired assertion. □

Note that the result in Lemma 2.1 in [2] (with $T_{\partial B(z,r)}$ instead of $\tau_{B(z,r)}$) may not be valid for our processes. This is because our processes might have killing inside $V$ and so $T_{\partial B(z,r)}$ may be different from $\tau_{B(z,r)}$.

LEMMA 2.9. *There exists $c > 0$ such that for every $r < r_1$ and $B(z,r) \subset \check{V}$,*

$$\mathbf{E}_z[\tau_{B(z,r)}] \vee \mathbf{E}_z[\widehat{\tau}_{B(z,r)}] < cr^2. \tag{2.18}$$

PROOF. By Proposition 2.7 and (2.6), the lemma is clear. In fact,

$$\mathbf{E}_z[\widehat{\tau}_{B(z,r)}] = \int_{B(z,r)} \frac{G_{B(z,r)}(y,z)H(y)}{H(z)}\,dy \le c\int_{B(z,r)} G^0_{B(z,r)}(z,y)\,dy \le c_1 r^2$$

for some constants $c, c_1 > 0$. □

Using the above lemma and the Markov property, we can easily get the following result.

LEMMA 2.10. *Suppose $r < r_1$, $B(z,r) \subset \check{V}$ and $U \subset D$. Then*

$$\mathbf{P}_z(\tau_U < \tau_{B(z,r)}) > c_1 \qquad (\textit{resp. } \mathbf{P}_z(\widehat{\tau}_U < \widehat{\tau}_{B(z,r)}) > c_1) \qquad \forall z$$

*for some $c_1 > 0$ implies*

$$\mathbf{E}_z[\tau_U] \le c_2 r^2 \qquad (\textit{resp. } \mathbf{E}_z[\widehat{\tau}_U] \le c_2 r^2) \qquad \forall z$$

*for some $c_2 > 0$.*

PROOF. Using (2.18) and the Markov property, the lemma can be proven using an argument similar to the one in the proof of Lemma 3.3 in [2] (with $\tau_U$ and $\tau_{B(z,r)}$ instead of the hitting times there). We omit the proof. □

**3. Parabolic and elliptic Harnack inequalities.** In this section, we shall prove a small-time parabolic Harnack inequality for $Y$ and $\widehat{Y}$. We will get a scale-invariant version of the elliptic Harnack inequality as a corollary. These Harnack inequalities will be used later to prove the main results of this paper.

Recall that $D$ is a bounded domain with $\overline{D} \subset V$, $\gamma_1 = \frac{1}{2}\mathrm{dist}(\partial V, \overline{D})$ and $\check{V} = \{z \in V; \rho_V(z) > \gamma_1\}$. In [17], we proved uniform Gaussian estimates for the density (with respect to the Lebesgue measure) of $Y^D$ when $D$ is a bounded smooth domain. We recall here part of the result from [17]: there exist positive constants $t_0$, $t_1$, $c_1$ and $c_2$ such that for every $R \le \sqrt{t_0}$, $t \le R^2 t_1$ and $(x, y) \in B(z, R) \times B(z, R)$,

$$\begin{aligned} &r^{B(z,R)}(t,x,y) \\ &\qquad \ge c_1 t^{-d/2}\left(1 \wedge \frac{\rho_{B(z,R)}(y)}{\sqrt{t}}\right)\left(1 \wedge \frac{\rho_{B(z,R)}(y)}{\sqrt{t}}\right)e^{-c_2|x-y|^2/(2t)}, \end{aligned} \tag{3.1}$$

whenever $B(z,R) \subset V$ (see Theorem 4.4(2) in [17]). In the remainder of this paper, $t_0$ and $t_1$ will always stand for the constants above.

With the density estimates (3.1) available, one can follow the ideas in [13] (see also [15, 27]) to prove the parabolic Harnack inequality. For this reason, the proofs of this section will be somewhat sketchy.

LEMMA 3.1. *For each $0 < \delta$, $u < 1$, there exists $\varepsilon = \varepsilon(d, \delta, u, t_1) > 0$ such that*

$$r^{B(x_0,R)}(t,x,y) \wedge \widehat{r}^{B(x_0,R)}(t,x,y) \ge \frac{\varepsilon}{|B(x_0,\delta R)|} \tag{3.2}$$

*for all $x, y \in B(x_0, \delta R) \subset \check{V}$, $R \le \sqrt{t_0}$ and $(1-u)R^2 t_1 \le t \le R^2 t_1$.*

PROOF. Fix $0 < \delta$, $u < 1$ and $B(x_0, \delta R) \subset \check{V}$. Let $B_R := B(x_0, R)$ and assume that $R \le \sqrt{t_0}$ and $t \le R^2 t_1$. By (2.6) and (3.1), there exist $c_1$ and $c_2$ such that

$$\begin{aligned} \widehat{r}^{B_R}(t,x,y) &= \frac{r^{B_R}(t,y,x)H(y)}{H(x)} \\ &\ge c_1 t^{-d/2}\left(1 \wedge \frac{\rho_{B_R}(y)}{\sqrt{t}}\right)\left(1 \wedge \frac{\rho_{B_R}(y)}{\sqrt{t}}\right)e^{-c_2|x-y|^2/(2t)}. \end{aligned} \tag{3.3}$$

If $|x - x_0| < \delta R$, $|y - x_0| < \delta R$ and $(1-u)R^2 t_1 \le t \le R^2 t_1$, then

$$\left(1 \wedge \frac{\rho_{B_R}(y)}{\sqrt{t}}\right)\left(1 \wedge \frac{\rho_{B_R}(y)}{\sqrt{t}}\right) \ge \frac{(1-\delta)^2}{t_1}$$

and

$$\frac{c_2|x-y|^2}{2t} \le \frac{2c_2\delta^2}{(1-u)t_1}.$$

So, the right-hand side of (3.3) is bounded below by

$$\begin{aligned} c_1(R^2t_1)^{-d/2}\frac{(1-\delta)^2}{t_1}e^{-2c_2\delta^2/((1-u)t_1)} &= c_3c_1t_1^{-d/2-1}\frac{(1-\delta)^2\delta^d}{|B(0,\delta R)|}e^{-2c_2\delta^2/((1-u)t_1)} \\ &=: \frac{\varepsilon}{|B(0,\delta)|}, \end{aligned}$$

where $c_3$ depends only on $d$. □

We define space-time processes $Z_s := (T_s, Y_s)$ and $\widehat{Z}_s := (T_s, \widehat{Y}_s)$, where $T_s = T_0 - s$. The law of the space-time processes $Z_s$ (and $\widehat{Z}_s$) starting from $(t,x)$ will be denoted by $\mathbf{P}_{t,x}$.

Definition 3.2. For any $(t,x) \in [0,\infty) \times V$, $u > 0$ and bounded subdomain $U$ of $V$, we say that a nonnegative continuous function $g$ defined on $[t, t+u] \times U$ is *parabolic* for $Y$ in $[t,t+u]\times U$ if for any $[s_1,s_2] \subset (t,t+u]$ and $B(y,\delta) \subset \overline{B(y,\delta)} \subset D$, we have

$$g(s,z) = \mathbf{E}_{s,z}[g(Z_{\tau_{(s_1,s_2]\times B(y,\delta)}}); Z_{\tau_{(s_1,s_2]\times B(y,\delta)}} \in (0,\infty)\times V] \tag{3.4}$$

for every $(s,z) \in (s_1,s_2] \times B(y,\delta)$, where $\tau_{(s_1,s_2]\times B(y,\delta)} = \inf\{s > 0 : Z_s \notin (s_1,s_2]\times B(y,\delta)\}$. The definition of parabolic functions for $\widehat{Y}$ is similar.

Lemma 3.3. *Suppose that $U$ is a subdomain of $V$. For each $T > 0$ and $y \in U$, $(t,x) \to r^U(t,x,y)$ and $(t,x) \to \widehat{r}^U(t,x,y)$ are parabolic in $(0,T]\times U$ for $Y$ and $\widehat{Y}$, respectively.*

Proof. See the proof of Lemma 4.5 in [6]. □

Corollary 3.4. *Suppose that $U$ is a subdomain of $V$. For each $T > 0$ and $y \in U$, and any nonnegative bounded function $f$ on $U$, the functions*

$$g(t,x) := \mathbf{E}_x[f(Y_t^U)] = \int_U r^U(t,x,y)f(y)\,dy$$

*and*

$$\widehat{g}(t,x) := \mathbf{E}_x[f(\widehat{Y}_t^U)] = \int_U \widehat{r}^U(t,x,y)f(y)\,dy$$

*are parabolic in $(0,T]\times U$ for $Y$ and $\widehat{Y}$, respectively.*

PROOF. The continuity of $\widehat{g}$ follows from the continuity of $\widehat{r}^U$. Equation (3.4) follows from Lemma 3.3 and Fubini's theorem. □

For $s \geq 0$, $R > 0$ and $B(x,R) \subset V$, we define the oscillation of a function $g$ on $(s - t_1R^2, s) \times B(x,R)$ by

$$\begin{aligned}&\mathrm{Osc}(g; s, x, R)\\ &\quad = \sup\{|g(s_1,x_1) - g(s_2,x_2)| : s_1, s_2 \in (s - t_1R^2, s),\ x_1, x_2 \in B(x,R)\}.\end{aligned}$$

LEMMA 3.5. *For any $0 < \delta < 1$, there exists $0 < \rho < 1$ such that for all $R \in (0, \sqrt{t_0}]$, $s \in [t_1R^2, \infty)$, $B(x_0,R) \subset \check{V}$ and function $g$ which is parabolic for $Y$ (resp. $\widehat{Y}$) in $(s - t_1R^2, s] \times B(x_0,R)$ and continuous in $[s - t_1R^2, s] \times \overline{B(x_0,R)}$,*

$$\mathrm{Osc}(g; s, x_0, \delta R) \leq \rho \mathrm{Osc}(g; s, x_0, R).$$

PROOF. Fix $s \geq 0$, $0 < R \leq \sqrt{t_0}$ and $B(x_0,R) \subset \check{V}$, and consider a function $g$ which is parabolic for $\widehat{Y}$ in $(s - t_1R^2, s] \times B(x_0,R)$ and continuous in $[s - t_1R^2, s] \times \overline{B(x_0,R)}$. Without loss of generality, we may assume that

$$\min_{(t,x)\in[s-t_1R^2,s]\times B(x_0,R)} g(t,x) = 0 \quad \text{and} \quad \max_{(t,x)\in[s-t_1R^2,s]\times B(x_0,R)} g(t,x) = 1.$$

Since $\widehat{Y}$ is a Hunt process, it is easy to see that $\widehat{Z}^\Omega$ is a Hunt process for any bounded open subset $\Omega$ of $[0,\infty) \times V$. So, $g$ and $1-g$ are excessive with respect to the process obtained by killing $\widehat{Z}$ upon exiting from $(s - t_1R^2, s) \times B(x_0,R)$. First, we assume that $\delta$ satisfies

$$\int_{B(x_0,\delta R)} g\left(s - \frac{1}{2}(\delta^2+1)t_1R^2, y\right) dy \geq \frac{|B(x_0,\delta R)|}{2}.$$

By Lemma 3.1, we have that for $(t,x) \in (s - \delta^2 t_1R^2, s) \times B(x_0,\delta R)$,

$$\begin{aligned} g(t,x) &\geq \mathbf{E}_{t,x}[g(\widehat{Z}_{t+1/2(\delta^2+1)t_1R^2-s}) : \\ &\qquad\qquad \widehat{Z}_{t+1/2(\delta^2+1)t_1R^2-s} \in (t_1R^2 - s, s) \times B(x_0,\delta R)] \\ &\geq \int_{B(x_0,\delta R)} \widehat{q}^{B(x_0,R)}\left(t + \frac{1}{2}(\delta^2+1)t_1R^2 - s, x, y\right) \\ &\qquad\qquad \times g\left(s - \frac{1}{2}(\delta^2+1)t_1R^2, y\right) dy \\ &\geq \frac{\varepsilon}{|B(x_0,\delta R)|}\frac{|B(x_0,\delta R)|}{2} = \frac{\varepsilon}{2}. \end{aligned}$$

Therefore, $\mathrm{Osc}(g; s, x_0, \delta R) \leq 1 - \varepsilon$.

If

$$\int_{B(x_0,\delta R)} g\left(s - \frac{1}{2}(\delta^2+1)t_1R^2, y\right) dy \leq \frac{|B(x_0, \delta R)|}{2},$$

then we consider $1-g$ and use the same argument as above. □

The above lemma implies the Hölder continuity of parabolic functions.

THEOREM 3.6. *For any $0<\delta<1$, there exist $c>0$ and $\beta \in (0,1)$ such that for all $R \in (0, \sqrt{t_0}]$, $s \in [t_1R^2, \infty)$, $B(x_0,R) \subset \check{V}$ and function $g$ which is parabolic for $Y$ (resp. $\widehat{Y}$) in $[s-t_1R^2, s] \times B(x_0,R)$ and continuous in $[s-t_1R^2, s] \times \overline{B(x_0,R)}$, we have*

$$|g(s_1,x_1) - g(s_2,x_2)| \leq c\|g\|_{L^\infty([s-t_1R^2,s]\times\overline{B(x_0,R)})}\left(\frac{|s_1-s_2|^2 + |x_1-x_2|}{R}\right)^\beta$$

*for any $(s_1,x_1),(s_2,x_2) \in [s - t_1\delta^2R^2, s] \times \overline{B(x_0,\delta R)}$.*

PROOF. See Theorem 5.3 in [13]. □

Using Lemmas 3.1 and 3.5, the proof of the next theorem is almost identical to that of Theorem 5.4 in [13]. Therefore, we omit the proof.

THEOREM 3.7. *For any $0<\alpha<\beta<1$ and $0<\delta<1$, there exists $c>0$ such that for all $R \in (0,\sqrt{t_0}]$, $s \in [t_1R^2,\infty)$, $B(x_0,R) \subset \check{V}$ positive function $\widehat{G}_D(x,y)$ is $\infty$ if and only if $X = y \in D$, and that for function $g$ which is parabolic for $Y$ (resp. $\widehat{Y}$) in $(s-t_1R^2, s] \times B(x_0,R)$ and continuous in $(s-t_1R^2, s] \times \overline{B(x_0,R)}$,*

$$g(t,y) \leq cg(s,x_0), \qquad (t,y) \in [s-\beta t_1R^2, s-\alpha t_1R^2] \times \overline{B(x_0,\delta R)}.$$

Now, the parabolic Harnack inequality is an easy corollary of the theorem above.

THEOREM 3.8 (Parabolic Harnack inequality). *For any $0<\alpha_1<\beta_1<\alpha_2<\beta_2<1$ and $0<\delta<1$, there exists $c>0$ such that for all $0<R\leq\sqrt{t_0}$, $B(x_0,R) \subset \check{V}$ and function $g$ which is parabolic for $Y$ (resp. $\widehat{Y}$) in $[0,t_1R^2) \times B(x_0,R)$ and continuous in $[0,t_1R^2] \times \overline{B(x_0,R)}$,*

$$\sup_{(t,y)\in B_1} g(t,y) \leq c \inf_{(t,y)\in B_2} g(t,y),$$

*where $B_i = \{(t,y) \in [\alpha_i t_1 R^2, \beta_i t_1 R^2] \times B(x_0,\delta R)\}$.*

The scale-invariant Harnack inequality is an easy corollary of the parabolic Harnack inequality.

THEOREM 3.9 (Scale-invariant Harnack inequality). *Every harmonic function for $Y$ (resp. $\widehat{Y}$) is Hölder continuous. There exists $c = c(D,V) > 0$ such that for every harmonic function $f$ for $Y$ (resp. $\widehat{Y}$) in $B(z_0,R)$ with $B(z_0,R) \subset \check{V}$, we have*

$$\sup_{y\in B(z_0,R/2)} f(y) \le c \inf_{y\in B(z_0,R/2)} f(y).$$

PROOF. By Proposition 2.6, any harmonic function $f$ for $\widehat{Y}$ in $B(z_0,r)$ is parabolic in $(0,T]\times B(z_0,r)$ with respect to $\widehat{Y}$ for any $T>0$. Thus, $f$ is Hölder continuous by Theorem 3.6 and the Harnack inequality above is true for small $R$ by Theorem 3.8. When $R$ is large and $B(z_0,R)\subset \check{V}$, we use a Harnack chain argument and the fact that $V$ is bounded. $\square$

**4. Analysis on various rough domains.** In this section, we recall the definitions of various rough domains from [1, 2] and prove the main lemma (Lemma 4.7). We will use the probabilistic methods used in [2]. For this reason, we follow the notation and the definitions of [2]. Unlike [2], we do not have the scaling property here and Lemma 2.8 works only for small balls. Moreover, our processes $Y$ and $\widehat{Y}$ may have killing inside $V$. All these make our argument more complicated than that of [2]. For the reader's convenience, we will spell out some of the proofs, especially the parts where things are more complicated.

A bounded domain $D$ is said to be a Hölder domain of order $\beta\in(0,1]$ if the boundary of $D$ is locally the graph of a function $\phi$ which is Hölder continuous of order $\beta$, that is, $|\phi(x)-\phi(z)|\le c|x-z|^{\beta}$. The concept of twisted Hölder domains, which is a natural generalization of the concept of Hölder domains, was introduced in [2]. Twisted Hölder domains have canals no longer and no thinner than Hölder domains, but do not have local representation of their boundaries as graphs of functions. For a rectifiable Jordan arc $\gamma$ and $x,y\in\gamma$, we denote the length of the piece of $\gamma$ between $x$ and $y$ by $l(\gamma(x,y))$. Recall the capacity defined in (2.12).

DEFINITION 4.1. A bounded domain $D\subset \mathbf{R}^d$ is called a *twisted Hölder domain of order* $\alpha\in(0,1]$ if there exist positive constants $c_1,\dots,c_5$, a point $z_0\in D$ and a continuous function $\delta: D\to(0,\infty)$ with the following properties:

(1) $\delta(x)\le \rho_D(x)^{\alpha}$ for all $x\in D$;

(2) for every $x \in D$, there exists a rectifiable Jordan arc $\gamma$ connecting $x$ and $z_0$ in $D$ such that

$$\delta(y) \geq c_2(l(\gamma(x,y)) + \delta(x)) \qquad \text{for all } y \in \gamma;$$

(3)

$$\frac{\mathrm{Cap}^{B(x,2c_3a)}(B(x,c_3a) \cap F(a)^c)}{\mathrm{Cap}^{B(x,2c_3a)}(B(x,c_3a))} \geq c_4 \qquad \text{for all } x \in F(a), a \leq c_5,$$

where $F(a) = \{y \in D : \delta(y) \leq a\}$.

One interesting fact is that the class of John domains (see page 422 of [2] for the definition) and the class of twisted Hölder domains of order 1 are identical (Proposition 3.2 of [2]). The boundary of a twisted Hölder domain can be highly nonrectifiable and, in general, no regularity of its boundary can be inferred. We refer to [2] for some elementary results on twisted Hölder domains.

Under some regularity assumption on the boundary of $D$, Bañuelos considered in [1] another natural generalization of Hölder domains. Let $k_D(x,y)$ be the quasi-hyperbolic distance

$$k_D(x,y) := \inf_{\gamma} \int_{\gamma} \frac{ds}{\rho_D(z)},$$

where the infimum is taken over all rectifiable curves joining $x$ to $y$ in $D$. The following definition is taken from [1].

Definition 4.2. A bounded domain $D \subset \mathbf{R}^d$ is called a *uniformly Hölder domain of order* $\alpha > 0$ if there exist positive constants $c_1, \ldots, c_5$ and a point $z_1 \in D$ with the following properties:

(1) $k_D(x, z_1) \leq c_1 \rho_D(x)^{-\alpha} + c_2$ for all $x \in D$;
(2) for every $Q \in \partial D$ and $r > 0$,

$$\mathrm{Cap}^{B(Q,2r)}(B(Q,r) \cap D^c) \geq c_3 r^{d-2}.$$

The class of uniformly Hölder domains is slightly more general than that of uniformly regular twisted $L^p$-domains defined in [2].

Lemma 4.3. (1) *If $D$ is a twisted Hölder domain of order $\alpha \in (0,1]$, there exist $c_1 > 0$, $a_1 > 0$ and $b_1 > 0$ such that for every $a \leq a_1$,*

$$\sup_{y \in F(a)} \mathbf{P}_y(T_{F(a)^c \cap B(y,ab_1)} < \tau_{B(y,2ab_1)}) > c_1$$

*and*

$$\sup_{y \in F(a)} \mathbf{P}_y(\widehat{T}_{F(a)^c \cap B(y,ab_1)} < \widehat{\tau}_{B(y,2ab_1)}) > c_1.$$

(2) *If $D$ is a uniformly Hölder domain of order $\alpha > 0$, there exist $c_2 > 0$ and $a_2 > 0$ such that for every $r \le a_2$,*

$$\sup_{y \in B(Q,2r/3) \cap D} \mathbf{P}_y(T_{B(Q,r)\cap D^c} < \tau_{B(Q,2r)}) > c_2$$

*and*

$$\sup_{y \in B(Q,2r/3) \cap D} \mathbf{P}_y(\widehat{T}_{B(Q,r)\cap D^c} < \widehat{\tau}_{B(Q,2r)}) > c_2.$$

Proof. Note that $\mathrm{Cap}^{B(x,2r)}(B(x,r)) \ge cr^{d-2}$. Thus, to prove (i), we only need to use Lemma 2.8 and Definition 4.1(3) with $K = \overline{B(0,1/2)}$ and $A^z = \partial F(a) \cap \overline{B(z,ab_1)} \subset \overline{B(z,ab_1)}$.

To prove (ii), we use Lemma 2.8 and Definition 4.2(2) with $K = \overline{B(0,2/3)}$ and $A^z = \overline{B(z,2r/3)} \cap \partial D \subset \overline{B(z,2r/3)}$. □

Definition 4.4. We say that a bounded domain $D \subset \mathbf{R}^d$ can be locally represented as the region above the graph of a function if there exist a positive constant $a_0$, a finite family of orthonormal coordinate systems $CS_j$'s, positive $b_j$'s and functions

$$f_j : \mathbf{R}^{d-1} \to (-\infty, 0], \quad j = 1, \dots, m_0,$$

such that

$$D = \bigcup_{j=1}^{m_0} \{x = (x_1, \dots, x_{d-1}, x_d) =: (\tilde{x}, x_d) \text{ in } CS_j : |\tilde{x}| < b_j, f_j(\tilde{x}) < x_d \le a_0\}.$$

Lemma 4.5. *Suppose that $D$ is a bounded domain which can be locally represented as the region above the graph of a function. Assume that $a \le r_1$ and that $y \in D$ is in $\{x = (\tilde{x}, x_d)$ in $CS_j : |\tilde{x}| < b_j, f_j(\tilde{x}) < x_d \le a_0\}$ for some $j = 1, \dots, m_0$. If $U$ and $M$ are subsets of $\mathbf{R}^d$ that can be written as*

$$U := \{(\tilde{x}, x_d) \text{ in } CS_j : |\tilde{x} - \tilde{y}| < a, |x_d - y_d| < a\},$$
$$M := \left\{(\tilde{x}, x_d) \text{ in } CS_j : |\tilde{x} - \tilde{y}| < \frac{a}{2}, x_d = a + y_d\right\},$$

*then there exists a constant $c_1 > 0$, independent of $a$, $y$ and $CS_j$, such that*

$$\left(\inf_{|\tilde{x}-\tilde{y}|<a/2, x_d = y_d} \mathbf{P}_x(T_M = \tau_U)\right) \wedge \left(\inf_{|\tilde{x}-\tilde{y}|<a/2, x_d = y_d} \mathbf{P}_x(\widehat{T}_M = \widehat{\tau}_U)\right) > c_1.$$

Proof. By our Harnack inequality (Theorem 3.9), it is enough to show that

$$\mathbf{P}_y(T_M \le \tau_U) \wedge \mathbf{P}_y(\widehat{T}_M \le \widehat{\tau}_U) > c_1$$

for some $c_1 > 0$ independent of $a$ and $CS_j$. Fix the coordinate systems $CS_j$. Let $B_1 := B(y, a)$ and

$$B_2 := B((\tilde{y}, y_d + a/2), a/\sqrt{2}),$$
$$M_1 := \left\{ (\tilde{x}, x_d); |\tilde{x} - \tilde{y}| < \frac{a}{2}, x_d = \frac{a}{2} + y_d \right\},$$
$$M_2 := \left\{ (\tilde{x}, x_d); |\tilde{x} - \tilde{y}| < \frac{a}{4}, x_d = a + y_d \right\}.$$

Note that $B_2 \cap \{x_d = a + y_d\} = M$. Thus,

$$\mathbf{P}_y(\widehat{T}_M = \widehat{\tau}_U) \geq \mathbf{E}_y[\mathbf{P}_{\widehat{Y}_{\widehat{T}_{M_1}}}(\widehat{T}_{M_2} < \widehat{\tau}_{B_2}); \widehat{T}_{M_1} < \widehat{\tau}_{B_1}]$$
$$\geq \mathbf{P}_y(\widehat{T}_{M_1} < \widehat{\tau}_{B_1}) \Big( \inf_{z \in M_1} \mathbf{P}_z(\widehat{T}_{M_2} < \widehat{\tau}_{B_2}) \Big).$$

Now, applying Lemma 2.8 to both factors on the right-hand side the equation above, we arrive at our desired conclusion. $\square$

For a bounded domain which can be locally represented as the region above the graph of a function, we put

$$\Theta := \frac{1}{2}\Big(1 + \frac{1}{4d-2}\Big).$$

For any $k < 0$ and $y \in D$ such that

$$y \in \{x = (\tilde{x}, x_d) \text{ in } CS_j : |\tilde{x}| < b_j, f_j(\tilde{x}) < x_d < 0\}$$

for some $j = 1, \ldots, m_0$, we let $l_0^{j,k}(y)$ be the smallest integer greater than $10|k|^{\Theta}(a_0/2 - y_d)/b_j$ and define

$$D_1^{j,k}(y) := \left\{ x \text{ in } CS_j : |\tilde{x} - \tilde{y}| < \frac{b_j}{4|k|^{\Theta}}, f_j(\tilde{x}) < x_d < a_0 \right\}, \tag{4.1}$$

$$D_2^{j,k}(y) := \left\{ x \text{ in } CS_j : |\tilde{x} - \tilde{y}| < \frac{b_j}{4|k|^{\Theta}}, |x_d - y_d| < \frac{b_j}{4|k|^{\Theta}} \right\}, \tag{4.2}$$

$$M^{j,k}(y) := \left\{ x \text{ in } CS_j : |\tilde{x} - \tilde{y}| < \frac{b_j}{20|k|^{\Theta}}, x_d = y_d + \frac{b_j}{10|k|^{\Theta}} l_0^{j,k}(y) \right\}, \tag{4.3}$$

where $a_0$, $b_j$, $CS_j$ and $f_j$ are the quantities from Definition 4.4.

Lemma 4.6. *Suppose that $D$ is a bounded domain which can be locally represented as the region above the graph of a function. There exists $p_0 \in (0, 1)$ such that if $p \in [p_0, 1)$ and*

$$k \leq - \max_{1 \leq j \leq m_0} \Big( \frac{b_j}{10 r_1} \Big)^{1/\Theta},$$

*then for any* $j = 1, \dots, m_0$,

$$\mathbf{P}_y(T_{\partial D} < \tau_{B(y,b_j|k|^{-\Theta})}) \le 1-p \qquad (\textit{resp. } \mathbf{P}_y(\widehat{T}_{\partial D} < \widehat{\tau}_{B(y,b_j|k|^{-\Theta})}) \le 1-p)$$

*for every* $y \in \{(\tilde{x}, x_d) \text{ in } CS_j : |\tilde{x}| < b_j, f_j(\tilde{x}) < x_d < 0\}$ *implies that*

$$\mathbf{P}_y(T_{M^{j,k}(y)} < \tau_{D_1^{j,k}(y)}) \ge \exp\Big(-c_1 \frac{8(a_0 - y_d)|k|^{\Theta}}{b_j}\Big)$$

$$\Big(\textit{resp. } \mathbf{P}_y(\widehat{T}_{M^{j,k}(y)} < \widehat{\tau}_{D_1^{j,k}(y)}) \ge \exp\Big(-c_1 \frac{8(a_0 - y_d)|k|^{\Theta}}{b_j}\Big)\Big)$$

*for some* $c_1 = c_1(p_0) > 0$ *independent of* $j, f_j$ *and* $y$.

Proof. Fix $j$ and $k$ satisfying the assumption of the lemma. We also fix a $y \in \{(\tilde{x}, x_d) \text{ in } CS_j : |\tilde{x}| < b_j, f_j(\tilde{x}) < x_d < 0\}$. Let $a := b_j 10^{-1}|k|^{-\Theta} \le r_1$ and

$$D_l := \Big\{x \text{ in } CS_j : |\tilde{x} - \tilde{y}| < \frac{5a}{2}, (y_d - a) \vee f_j(\tilde{x}) < x_d < y_d + al\Big\}, \qquad l \ge 1,$$

$$\tilde{D}_l := \Big\{x \text{ in } CS_j : |\tilde{x} - \tilde{y}| < \frac{5a}{2}, y_d - a < x_d < y_d + al\Big\}, \qquad l \ge 1,$$

$$W_l := \Big\{x \text{ in } CS_j : |\tilde{x} - \tilde{y}| < \frac{5a}{2}, y_d + a(l-5) \le x_d < y_d + al\Big\}, \qquad l \ge 4,$$

$$V_l := \{x \text{ in } CS_j : x \in D_l^c, |\tilde{x} - \tilde{y}| < a, y_d + a(l-2) < x_d < y_d + al\}, \qquad l \ge 3,$$

$$B_l := \Big\{x \text{ in } CS_j ; \Big|x - \Big(\tilde{y}, a\Big(l - \frac{5}{2}\Big)\Big)\Big| < \frac{5a}{2}\Big\}, \qquad l \ge 4$$

and

$$2B_l := \Big\{x \text{ in } CS_j ; \Big|x - \Big(\tilde{y}, a\Big(l - \frac{5}{2}\Big)\Big)\Big| < 5a\Big\}, \qquad l \ge 4.$$

Note that $V_l \subset B_{l+1} \subset W_{l+1} \subset 2B_{l+1}$ (see Figure 1). Since $\mathbf{P}_y(\widehat{T}_{\partial D} < \widehat{\tau}_{B(y,b_j|k|^{-\Theta})}) \le 1-p$, we have

$$\mathbf{P}_y(\widehat{T}_{V_3} < \widehat{\tau}_{B_4}) \le \mathbf{P}_y(\widehat{T}_{\partial D_4} < \widehat{\tau}_{\tilde{D}_4}) \le (1-p). \tag{4.4}$$

Thus, by Lemma 2.8,

$$\mathrm{Cap}^{B_4}(V_3) \le c_1 \Big(\inf_{w \in K} G^0_{B(0,1)}((0, -\tfrac{3}{5}), w)\Big)^{-1} a^{d-2}(1-p) \le c_2 a^{d-2}(1-p)$$

for some constants $c_1, c_2 > 0$ and where

$$K := \{|\tilde{x}| < 1/10, -1/5 < x_d < 3/5\}. \tag{4.5}$$

By the translation invariance of Cap and the definition of $V_l$,

$$\mathrm{Cap}^{B_{l+1}}(V_l) \le c_2 a^{d-2}(1-p).$$

Since $W_{l+1} \subset 2B_{l+1}$, by Lemma 2.8, for $y_l := y + (\tilde{0}, (l-3)a)$,

$$\begin{aligned}\mathbf{P}_{y_l}(\widehat{T}_{V_l} < \widehat{\tau}_{W_{l+1}}) &\le \mathbf{P}_{y_l}(\widehat{T}_{V_l} < \tau_{2B_{l+1}}) \\ &\le c_3 a^{2-d}\Big(\sup_{w \in 1/2K} G^0_{B(0,1)}((0, -\tfrac{3}{10}), w)\Big)\,\mathrm{Cap}^{2B_{l+1}}(V_l) \\ &\le c_4 a^{2-d}\,\mathrm{Cap}^{2B_{l+1}}(V_l),\end{aligned}$$

where $K$ is defined in (4.5). But, by the definition of Cap, $\mathrm{Cap}^{2B_{l+1}}(V_l) \le \mathrm{Cap}^{B_{l+1}}(V_l)$. Therefore,

$$\mathbf{P}_{y_l}(\widehat{T}_{V_l} < \widehat{\tau}_{W_{l+1}}) \le c_4 a^{2-d}\,\mathrm{Cap}^{B_{l+1}}(V_l) \le c_5(1-p).$$

Applying the Harnack inequality (Theorem 3.9), we get

$$\mathbf{P}_x(\widehat{T}_{V_l} < \widehat{\tau}_{W_{l+1}}) \le c_6(1-p), \qquad |\tilde{x} - \tilde{y}| < \frac{a}{2}, \qquad x_d = y_d + a(l-3). \tag{4.6}$$

Using our Lemma 4.5 and (4.6) instead of Lemma 2.3 and (2.5) of [2], the remaining part of the proof is similar to the proof of Lemma 2.4 on page 414, starting from the line 3, in [2] (after rescaling) with

$$\hat{D}_l := \{x \text{ in } CS_j : |\tilde{x} - \tilde{y}| < a, y_d - a < x_d < y_d + al\},$$

$$M_l := \Big\{x \text{ in } CS_j : |\tilde{x} - \tilde{y}| < \frac{a}{2}, x_d = y_d + al\Big\}.$$

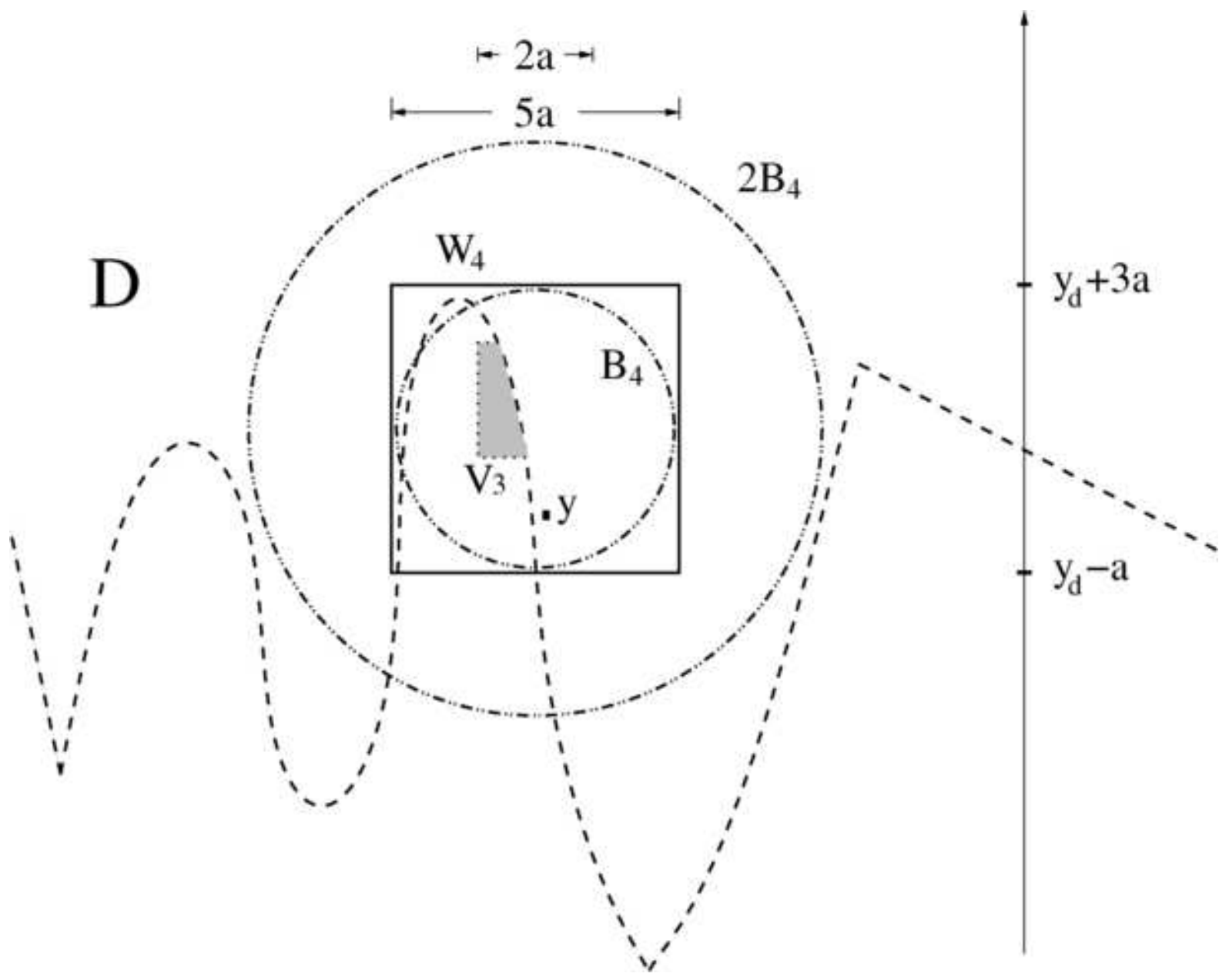

Fig. 1. $V_3 \subset B_4 \subset W_4 \subset 2B_4$.

However, due to the possible killing inside the domain in our case, things are more delicate. We include the details of the remaining part of the proof for the reader's convenience.

Let $\theta$ be the usual shift operator for Markov processes and define

$$A_l := \bigcap_{m=1}^{l} \{\widehat{\tau}_{D_m} = \widehat{T}_{M_m}, \widehat{T}_{\partial \hat{D}_{m-1}} \circ \theta_{\widehat{T}_{M_m}} > \widehat{\tau}_{\hat{D}_l}\}.$$

Note that by the strong Markov property applied at $\widehat{T}_{M_1}$,

$$\begin{aligned}
&\mathbf{P}_y\left(\bigcap_{m=1}^{4} \{\widehat{\tau}_{\hat{D}_m} = \widehat{T}_{M_m}, \widehat{T}_{\partial \hat{D}_{m-1}} \circ \theta_{\widehat{T}_{M_m}} > \widehat{\tau}_{\hat{D}_4}\}\right) \\
&\quad = \mathbf{E}_y\Bigg[\mathbf{P}_{\widehat{Y}^D_{\widehat{T}_{M_1}}}\Bigg(\bigcap_{m=2}^{4} \{\widehat{\tau}_{\hat{D}_m} = \widehat{T}_{M_m}, \widehat{T}_{\partial \hat{D}_{m-1}} \circ \theta_{\widehat{T}_{M_m}} > \widehat{\tau}_{\hat{D}_4}\} \\
&\qquad\qquad \cap \{\widehat{T}_{\partial \hat{D}_0} > \widehat{\tau}_{\hat{D}_4}\}\Bigg) : \widehat{\tau}_{\hat{D}_1} = \widehat{T}_{M_1}\Bigg].
\end{aligned}$$

Thus, by Lemma 4.5 and the strong Markov property applied at $\widehat{T}_{M_m}$, $m = 1, \dots, 4$, we get

$$\begin{aligned}
&\mathbf{P}_y\left(\bigcap_{m=1}^{4} \{\widehat{\tau}_{\hat{D}_m} = \widehat{T}_{M_m}, \widehat{T}_{\partial \hat{D}_{m-1}} \circ \theta_{\widehat{T}_{M_m}} > \widehat{\tau}_{\hat{D}_4}\}\right) \\
&\quad \geq c_7 \inf_{x \in M_1} \mathbf{P}_x\left(\bigcap_{m=2}^{4} \{\widehat{\tau}_{\hat{D}_m} = \widehat{T}_{M_m}, \widehat{T}_{\partial \hat{D}_{m-1}} \circ \theta_{\widehat{T}_{M_m}} > \widehat{\tau}_{\hat{D}_4}\} \cap \{\widehat{T}_{\partial \hat{D}_0} > \widehat{\tau}_{\hat{D}_4}\}\right) \\
&\quad = c_7 \inf_{x \in M_1} \mathbf{E}_x\Bigg[\mathbf{P}_{\widehat{Y}^D_{\widehat{\tau}_{M_2}}}\Bigg(\bigcap_{m=3}^{4} \{\widehat{\tau}_{\hat{D}_m} = \widehat{T}_{M_m}, \\
&\qquad\qquad \widehat{T}_{\partial \hat{D}_{m-1}} \circ \theta_{\widehat{T}_{M_m}} > \widehat{\tau}_{\hat{D}_4}\} \qquad (4.7)\\
&\qquad\qquad \cap \{\widehat{T}_{\partial \hat{D}_0} > \widehat{\tau}_{\hat{D}_4}, \widehat{T}_{\partial \hat{D}_1} > \widehat{\tau}_{\hat{D}_4}\}\Bigg) : \widehat{\tau}_{\hat{D}_2} = \widehat{T}_{M_2}\Bigg] \\
&\quad \geq c_7^2 \cdots \\
&\quad \geq c_7^4 \inf_{x \in M_4} \mathbf{P}_x\left(\bigcap_{m=1}^{4} \{\widehat{T}_{\partial \hat{D}_{m-1}} > \widehat{\tau}_{\hat{D}_4}\}\right) \\
&\quad = c_7^4 \inf_{x \in M_4} \mathbf{P}_x\left(\bigcap_{m=1}^{4} \{\widehat{T}_{\partial \hat{D}_{m-1}} > 0\}\right) = c_7^4.
\end{aligned}$$

On the other hand, since

$$\{\widehat{\tau}_{\tilde{D}_4} \le \widehat{T}_{\partial D_4}, \widehat{\tau}_{\hat{D}_4} = \widehat{T}_{M_4}\} = \{\widehat{\tau}_{\tilde{D}_4} \le \widehat{T}_{\partial D_4} \le \widehat{T}_{M_4} = \widehat{\tau}_{\hat{D}_4} \le \widehat{\tau}_{\tilde{D}_4}\}$$
$$= \{\widehat{\tau}_{\tilde{D}_4} = \widehat{T}_{\partial D_4} = \widehat{T}_{M_4}\} \subset \{\widehat{\tau}_{D_4} = \widehat{T}_{M_4}\},$$

we have

$$\begin{aligned}
&\mathbf{P}_y\left(\bigcap_{m=1}^{4}\{\widehat{\tau}_{\hat{D}_m} = \widehat{T}_{M_m}, \widehat{T}_{\partial\hat{D}_{m-1}} \circ \theta_{\widehat{T}_{M_m}} > \widehat{\tau}_{\hat{D}_4}\}\right) \\
&\quad = \mathbf{P}_y\left(\bigcap_{m=1}^{4}\{\widehat{\tau}_{\hat{D}_m} = \widehat{T}_{M_m}, \widehat{T}_{\partial\hat{D}_{m-1}} \circ \theta_{\widehat{T}_{M_m}} > \widehat{\tau}_{\hat{D}_4}\} \cap \{\widehat{\tau}_{\tilde{D}_4} \le \widehat{T}_{\partial D_4}\}\right) \\
&\qquad + \mathbf{P}_y\left(\bigcap_{m=1}^{4}\{\widehat{\tau}_{\hat{D}_m} = \widehat{T}_{M_m}, \widehat{T}_{\partial\hat{D}_{m-1}} \circ \theta_{\widehat{T}_{M_m}} > \widehat{\tau}_{\hat{D}_4}\} \cap \{\widehat{\tau}_{\tilde{D}_4} > \widehat{T}_{\partial D_4}\}\right) \\
&\quad \le \mathbf{P}_y\left(\bigcap_{m=1}^{3}\{\widehat{\tau}_{\hat{D}_m} = \widehat{T}_{M_m}, \widehat{T}_{\partial\hat{D}_{m-1}} \circ \theta_{\widehat{T}_{M_m}} > \widehat{\tau}_{\hat{D}_4}\} \cap \{\widehat{\tau}_{D_4} = \widehat{T}_{M_4}\}\right) \\
&\qquad + \mathbf{P}_y(\widehat{\tau}_{\tilde{D}_4} > \widehat{T}_{\partial D_4}) \\
&\quad \le \mathbf{P}_y\left(\bigcap_{m=1}^{3}\{\widehat{\tau}_{D_m} = \widehat{T}_{M_m}, \widehat{T}_{\partial\hat{D}_{m-1}} \circ \theta_{\widehat{T}_{M_m}} > \widehat{\tau}_{\hat{D}_4}\} \cap \{\widehat{\tau}_{D_4} = \widehat{T}_{M_4}\}\right) \\
&\qquad + \mathbf{P}_y(\widehat{\tau}_{\tilde{D}_4} > \widehat{T}_{\partial D_4}) \\
&\quad = \mathbf{P}_y(A_4) + \mathbf{P}_y(\widehat{\tau}_{\tilde{D}_4} > \widehat{T}_{\partial D_4}) \le \mathbf{P}_y(A_4) + 1 - p.
\end{aligned}$$

In the last inequality above, we have used (4.4). Letting $p > 1 - c_7^4/2$ and combining the inequality above with (4.7), we have

$$\mathbf{P}_y(A_2) \ge \mathbf{P}_y(A_3) \ge \mathbf{P}_y(A_4) \ge c_7^4/2. \tag{4.8}$$

We claim that there exist $c_8$ and $p_0$, which will be chosen later, such that for every $p > p_0$,

$$\mathbf{P}_y(A_{l+1}) \ge c_8 \mathbf{P}_y(A_l), \qquad l \ge 2. \tag{4.9}$$

We will prove this claim by induction. By (4.8), we know that the claim is valid for $l = 2, 3$. First, we note that by Lemma 4.5 and the strong Markov property applied at $\widehat{T}_{M_{l+1}}$, we get

$$\begin{aligned}
&\mathbf{P}_y(A_{l+1} \cap \{\widehat{\tau}_{\hat{D}_{l+2}} = \widehat{T}_{M_{l+2}}, \widehat{T}_{\partial\hat{D}_l} \circ \theta_{\widehat{T}_{M_{l+1}}} > \widehat{\tau}_{\hat{D}_{l+2}}\}) \\
&\quad = \mathbf{E}_y[\mathbf{P}_{\widehat{Y}^D_{\widehat{T}_{M_{l+1}}}}(\widehat{\tau}_{\hat{D}_{l+2}} = \widehat{T}_{M_{l+2}}, \widehat{T}_{\partial\hat{D}_l} > \widehat{\tau}_{\hat{D}_{l+2}}) : A_{l+1}]
\end{aligned}$$

$$\geq \inf_{x\in M_{l+1}} \mathbf{P}_x(\widehat{\tau}_{\hat{D}_{l+2}} = \widehat{T}_{M_{l+2}}, \widehat{T}_{\partial\hat{D}_l} > \widehat{\tau}_{\hat{D}_{l+2}})\mathbf{P}_y(A_{l+1})$$
$$\geq c_7\mathbf{P}_y(A_{l+1}).$$

On the other hand,

$$\begin{aligned}
&\mathbf{P}_y(A_{l+1}\cap\{\widehat{\tau}_{\hat{D}_{l+2}} = \widehat{T}_{M_{l+2}}, \widehat{T}_{\partial\hat{D}_l}\circ\theta_{\widehat{T}_{M_{l+1}}} > \widehat{\tau}_{\hat{D}_{l+2}}\})\\
&\quad= \mathbf{P}_y(A_{l+1}\cap\{\widehat{\tau}_{D_{l+2}} = \widehat{T}_{M_{l+2}} = \widehat{\tau}_{\hat{D}_{l+2}}, \widehat{T}_{\partial\hat{D}_l}\circ\theta_{\widehat{T}_{M_{l+1}}} > \widehat{\tau}_{\hat{D}_{l+2}}\})\\
&\qquad+\mathbf{P}_y(A_{l+1}\cap\{\widehat{\tau}_{D_{l+2}} \neq \widehat{T}_{M_{l+2}}, \widehat{\tau}_{\hat{D}_{l+2}} = \widehat{T}_{M_{l+2}}, \widehat{T}_{\partial\hat{D}_l}\circ\theta_{\widehat{T}_{M_{l+1}}} > \widehat{\tau}_{\hat{D}_{l+2}}\})\\
&\quad= \mathbf{P}_y\Bigg(\bigcap_{m=1}^{l+1}\{\widehat{\tau}_{D_m} = \widehat{T}_{M_m}, \widehat{T}_{\partial\hat{D}_{m-1}}\circ\theta_{\widehat{T}_{M_m}} > \widehat{\tau}_{\hat{D}_{l+1}}\}\\
&\qquad\qquad\cap\{\widehat{\tau}_{D_{l+2}} = \widehat{T}_{M_{l+2}} = \widehat{\tau}_{\hat{D}_{l+2}}, \widehat{T}_{\partial\hat{D}_l}\circ\theta_{\widehat{T}_{M_{l+1}}} > \widehat{\tau}_{\hat{D}_{l+2}}\}\Bigg)\\
&\qquad+\mathbf{P}_y\Bigg(\bigcap_{m=1}^{l-1}\{\widehat{\tau}_{D_m} = \widehat{T}_{M_m}, \widehat{T}_{\partial\hat{D}_{m-1}}\circ\theta_{\widehat{T}_{M_m}} > \widehat{\tau}_{\hat{D}_{l+1}}\}\\
&\qquad\qquad\cap\{\widehat{\tau}_{D_l} = \widehat{T}_{M_l}, \widehat{T}_{\partial\hat{D}_{l-1}}\circ\theta_{\widehat{T}_{M_l}} > \widehat{\tau}_{\hat{D}_{l+1}}, \widehat{\tau}_{D_{l+1}} = \widehat{T}_{M_{l+1}},\\
&\qquad\qquad\quad\widehat{T}_{\partial\hat{D}_l}\circ\theta_{\widehat{T}_{M_{l+1}}} > \widehat{\tau}_{\hat{D}_{l+2}}, \widehat{\tau}_{D_{l+2}} \neq \widehat{T}_{M_{l+2}}, \widehat{\tau}_{\hat{D}_{l+2}} = \widehat{T}_{M_{l+2}}\}\Bigg)\\
&\quad\leq \mathbf{P}_y\Bigg(\bigcap_{m=1}^{l+1}\{\widehat{\tau}_{D_m} = \widehat{T}_{M_m}, \widehat{T}_{\partial\hat{D}_{m-1}}\circ\theta_{\widehat{T}_{M_m}} > \widehat{\tau}_{\hat{D}_{l+2}}\}\\
&\qquad\qquad\cap\{\widehat{\tau}_{D_{l+2}} = \widehat{T}_{M_{l+2}} = \widehat{\tau}_{\hat{D}_{l+2}}\}\Bigg)\\
&\qquad+\mathbf{P}_y\Bigg(\bigcap_{m=1}^{l-1}\{\widehat{\tau}_{D_m} = \widehat{T}_{M_m}, \widehat{T}_{\partial\hat{D}_{m-1}}\circ\theta_{\widehat{T}_{M_m}} > \widehat{\tau}_{\hat{D}_{l-1}}\}\\
&\qquad\qquad\cap\{\widehat{\tau}_{D_l} = \widehat{T}_{M_l}, \widehat{T}_{\partial\hat{D}_{l-1}}\circ\theta_{\widehat{T}_{M_l}} > \widehat{\tau}_{\hat{D}_{l+1}}, \widehat{\tau}_{D_{l+1}} = \widehat{T}_{M_{l+1}},\\
&\qquad\qquad\quad\widehat{T}_{\partial\hat{D}_l}\circ\theta_{\widehat{T}_{M_{l+1}}} > \widehat{\tau}_{\hat{D}_{l+2}}, \widehat{\tau}_{D_{l+2}} < \widehat{\tau}_{\hat{D}_{l+2}} = \widehat{T}_{M_{l+2}}\}\Bigg)\\
&\quad\leq \mathbf{P}_y(A_{l+2}) + \mathbf{P}_y(A_{l-1}\cap\{\widehat{T}_{V_{l+2}}\circ\theta_{\widehat{T}_{M_{l-1}}} < \widehat{\tau}_{W_{l+3}}\circ\theta_{\widehat{T}_{M_{l-1}}}\}),
\end{aligned}$$

which is less than or equal to $\mathbf{P}_y(A_{l+2}) + c_6(1-p)\mathbf{P}_y(A_{l-1})$, by (4.6). Combining the two inequalities above, we get, by induction,

$$\begin{aligned}\mathbf{P}_y(A_{l+2}) &\geq c_7\mathbf{P}_y(A_{l+1}) - c_6(1-p)\mathbf{P}_y(A_{l-1}) \\ &\geq c_7\mathbf{P}_y(A_{l+1}) - c_6(1-p)c_8^{-2}\mathbf{P}_y(A_{l+1}) \\ &= (c_7 - c_6(1-p)c_8^{-2}2)\mathbf{P}_y(A_{l+1}).\end{aligned}$$

Choose $c_8 < c_7^4/2$ small and then choose $p_0 < 1$ large so that for every $p \in [p_0, 1)$,

$$c_7 - c_6(1-p)c_8^2 > c_8.$$

Thus, the claim (4.9) is valid. Recall that $l_0 := l_0^{j,k}(y)$ is the smallest integer greater than $(a_0/2 - y_d)/a$. From (4.8) and (4.9), we conclude that

$$\begin{aligned}\mathbf{P}_y(\widehat{T}_{M^{j,k}(y)} < \widehat{\tau}_{D_1^{j,k}(y)}) &\geq \mathbf{P}_y(\widehat{\tau}_{D_{l_0}} = \widehat{T}_{M_{l_0}}) \geq \mathbf{P}_y(A_{l_0}) \\ &\geq c_8^{l_0-2}\mathbf{P}_y(A_2) \geq \frac{c_7^4}{2}c_8^{l_0-2} \geq \exp\left(-c_9\frac{8(a_0-y_d)|k|^\Theta}{b_j}\right)\end{aligned}$$

for some positive constant $c_9$. $\square$

For any positive function $h$ which is harmonic in $D$ for either $Y$ or $\widehat{Y}$, we let $S_k := \{x \in D : h(x) \leq 2^{k+1}\}$.

Lemma 4.7. *Suppose that $D$ is one of the following types of bounded domains:*

(a) *a twisted Hölder domain of order $\alpha \in (1/3, 1]$;*
(b) *a uniformly Hölder domain of order $\alpha \in (0, 2)$;*
(c) *a domain which can be locally represented as the region above the graph of a function.*

*Then, for any positive bounded function $h$ which is harmonic in $D$ for $Y$ (resp. $\widehat{Y}$), there exist $c > 0$ and $\beta > 0$ such that*

$$\sup_{x\in D}\mathbf{E}_x[\tau_{S_k}] \leq c|k|^{-1-\beta} \qquad \left(\textit{resp. } \sup_{x\in D}\mathbf{E}_x[\widehat{\tau}_{S_k}] \leq c|k|^{-1-\beta}\right). \tag{4.10}$$

Proof. Note that by (2.4) and (2.6), we have

$$\widehat{G}_D(x,y) = \frac{G_D(y,x)H(y)}{H(x)} \leq c|x-y|^{-d+2},$$

which implies that

$$\sup_{x\in D}\mathbf{E}_x[\widehat{\tau}_{S_k}] \leq \sup_{x\in D}\mathbf{E}_x[\widehat{\tau}_D] \leq c_1 \sup_{x\in D}\int_D |x-y|^{-d+2}\,dy < \infty.$$

Thus, we only need to show (4.10) for negative $k$ with $|k|$ large.

(i) Assume $D$ is a twisted Hölder domain of order $\alpha \in (1/3, 1)$. Recall that $z_0$ is the point from Definition 4.1(2). By Lemma 3.1 in [2], there exists $c_1 = c_1(D) > 0$ such that for every $x \in D$, there exists a sequence of open balls contained in $D$, with centers $z^1 = x, z^2, \dots, z^k = z_0$ and radii $a_j \le \text{dist}(z^j, \partial D)$, such that $|z^j - z^{j+1}| < (a_j \wedge a_{j+1})/2$ and $k \le c_1 \delta(x)^{1-1/\alpha}$. Thus, by the Harnack inequality (Theorem 3.9), there exists $c_2 = c_2(z_0) > 0$ such that

$$h(x) \ge \exp(-c_2 \delta(x)^{1-1/\alpha}). \tag{4.11}$$

If $x \in S_k$, then from (4.11), we have

$$2^{k+1} \ge h(x) \ge \exp(-c_2 \delta(x)^{1-1/\alpha}),$$

which implies that there exists $c_3 > 0$ such that

$$\delta(x) \le c_3 |k|^{-\alpha/(1-\alpha)}.$$

Therefore, $S_k \subset F(a)$ with $a \le c_3 |k|^{-\alpha/(1-\alpha)}$. We consider negative $k$ with $|k|$ large enough so that

$$c_3 |k|^{-\alpha/(1-\alpha)} \le a_1 \quad \text{and} \quad 2c_3 b_1 |k|^{-\alpha/(1-\alpha)} \le |k|^{-\alpha+1/(4(1-\alpha))},$$

where $a_1$ and $b_1$ are the constants in Lemma 4.3. Note that the above is always possible because $\frac{1}{4}(\alpha+1) < \alpha$. For those $k$, we apply Lemma 4.3 and get

$$\begin{aligned}
&\mathbf{P}_x(\widehat{\tau}_{S_k} < \widehat{\tau}_{B(x,|k|^{-(\alpha+1)/(4(1-\alpha))})}) \\
&\quad \ge \mathbf{P}_x(\widehat{\tau}_{S_k} < \widehat{\tau}_{B(x,2c_3b_1|k|^{-\alpha/(1-\alpha)})}) \\
&\quad \ge \mathbf{P}_x(\widehat{\tau}_{F(c_3|k|^{-\alpha/(1-\alpha)})} < \widehat{\tau}_{B(x,2c_3b_1|k|^{-\alpha/(1-\alpha)})}) \\
&\quad \ge \mathbf{P}_x(\widehat{T}_{F(c_3|k|^{-\alpha/(1-\alpha)})^c \cap B(x,c_3b_1|k|^{-\alpha/(1-\alpha)})} < \widehat{\tau}_{B(x,2c_3b_1|k|^{-\alpha/(1-\alpha)})}) \ge c_4
\end{aligned}$$

for some $c_4 > 0$. Thus, by Lemma 2.10, we have

$$\mathbf{E}_x[\widehat{\tau}_{S_k}] \le c_5 |k|^{-\alpha+1/(2(1-\alpha))} = c_5 |k|^{-1-\beta},$$

where $\beta = (3\alpha - 1)/(2 - 2\alpha) > 0$.

(ii) Assume that $D$ is a John domain (i.e., a twisted Hölder domain of order $\alpha = 1$). It is well known that $k_D(x, z_0) \le -c_6 \ln \rho_D(x) + c_7$ for some positive constants $c_6, c_7$ (see, e.g., page 185 in [1]). It is easy to see that the shortest length of a Harnack chain connecting $x$ and $z_1$ is comparable to $k_D(x, z_0)$. Thus, by our Harnack inequality (Theorem 3.9),

$$h(x) \ge \exp(-c_8 k_D(x, z_1)) \ge c_9 \rho_D(x)^{c_{10}} \ge c_{11} \delta(x)^{c_{10}}$$

for some positive constants $c_8, c_9, c_{10}$. Using the above instead of (4.11), we can repeat the argument in (i) to arrive at the desired conclusion. We omit the details.

(iii) We now assume that $D$ is a uniformly Hölder domain of order $\alpha \in (0,2)$. Recall that $z_1$ is the point from Definition 4.2. Since the shortest length of a Harnack chain connecting $x$ and $z_1$ is comparable to $k_D(x,z_1)$, by the Harnack inequality (Theorem 3.9) and Definition 4.2(1), there exists $c_{11} = c_{11}(z_1) > 0$ such that

$$h(x) \geq \exp(-ck_D(x,z_1)) \geq \exp(-c_{11}\rho_D(x)^{-\alpha}). \tag{4.12}$$

If $x \in S_k$, from (4.12), we have

$$2^{k+1} \geq h(x) \geq \exp(-c_{11}\rho_D(x)^{-\alpha}),$$

which implies that there exists $c_{12} > 0$ such that

$$\rho_D(x) \leq c_{12}|k|^{-1/\alpha}.$$

Therefore, $S_k \subset D(a) := \{x \in D : \rho_D(x) < a\}$ with $a \leq c_{12}|k|^{-1/\alpha}$. For each $x \in S_k$, choose a point $Q_x \in \partial D$ such that

$$|Q_x - x| = \frac{3c_{12}}{2}|k|^{-1/\alpha}.$$

We consider negative $k$ with $|k|$ large enough such that

$$c_{12}|k|^{-1/\alpha} \leq a_2 \quad \text{and} \quad \tfrac{7}{2}c_{12}|k|^{-1/\alpha} \leq |k|^{-(\alpha+2)/(4\alpha)},$$

where $a_2$ is the constant in Lemma 4.3(2). Note that the above is always possible because $\frac{1}{4}(\alpha+2) < 1$. We also note that for those negative $k$'s,

$$B(Q_x, 2c_{12}|k|^{-1/\alpha}) \subset B(x, |k|^{-(\alpha+2)/(4\alpha)}).$$

For those negative $k$'s, we apply Lemma 4.3 and get

$$\begin{aligned}
&\mathbf{P}_x(\widehat{\tau}_{S_k} < \widehat{\tau}_{B(x,|k|^{-(\alpha+2)/(4\alpha)})}) \\
&\quad \geq \mathbf{P}_x(\widehat{\tau}_{S_k} < \widehat{\tau}_{B(Q_x,2c_{12}|k|^{-1/\alpha})}) \\
&\quad \geq \mathbf{P}_x(\widehat{\tau}_{D(c_{12}|k|^{-1/\alpha})} < \widehat{\tau}_{B(Q_x,2c_{12}|k|^{-1/\alpha})}) \\
&\quad = \mathbf{P}_x(\widehat{\tau}_{D(c_{12}|k|^{-1/\alpha})\cap B(Q_x,2c_{12}|k|^{-1/\alpha})} < \widehat{\tau}_{B(Q_x,2c_{12}|k|^{-1/\alpha})}) \\
&\quad \geq \mathbf{P}_x(\widehat{T}_{D^c\cap B(Q_x,c_{12}|k|^{-1/\alpha})} < \widehat{\tau}_{B(Q_x,2c_{12}|k|^{-1/\alpha})}) \geq c_4
\end{aligned}$$

for some constant $c_{13} > 0$. Thus, by Lemma 2.10, we have

$$\mathbf{E}_x[\widehat{\tau}_{S_k}] \leq c_{14}|k|^{-(\alpha+2)/(2\alpha)} = c_{14}|k|^{-1-\beta}$$

for some constant $c_{14} > 0$, where $\beta = \frac{1}{2}(2-\alpha)/\alpha > 0$.

(iv) Finally, we assume that $D$ is a bounded domain which can be locally represented as the region above the graph of a function. Without loss of generality, we may assume that $\max_{1\le i\le m_0} f_i < -\varepsilon$ for some positive $\varepsilon = \varepsilon(D)$ so that

$$\bigcup_{i=1}^{m_0}\{(\tilde{x},x_d) \text{ in } CS_i : |\tilde{x}|<b_i, 0\le x_d\le a_0\}$$

is a compact subset of $D$. Thus, by the continuity of $h$, there exists $k_0>0$ such that $h(x)\ge 2^{-k_0+1}$ for $x\in K$. We let

$$k_1 := k_0 \vee \max_{1\le i\le m_0}\left(\frac{b_i}{r_1}\right)^{1/\Theta} \qquad \text{where } \Theta = \frac{1}{2}\left(1+\frac{1}{4d-2}\right).$$

Fix $j$ and $f_j$, and consider $y\in\{(\tilde{x},x_d) \text{ in } CS_j : |\tilde{x}|<b_j, f_j(\tilde{x})<x_d<0\}$. Recall that $p_0$ is the constant in Lemma 4.6.

We claim that there exist $p_1\in(p_0,1)$ and $k_2\ge k_1$ such that for every

$$y\in S_k\cap\{(\tilde{x},x_d) \text{ in } CS_j : |\tilde{x}|<b_j, f_j(\tilde{x})<x_d<0\}$$

and $k<-k_2$, we have

$$\mathbf{P}_y(\widehat{T}_{\partial D}<\widehat{\tau}_{B(y,b_j|k|^{-\Theta})})>1-p_1.$$

Recall that $D_1^{j,k}(y)$, $D_2^{j,k}(y)$ and $M^{j,k}(y)$ are defined in (4.1)–(4.3). If we suppose that

$$\mathbf{P}_y(\widehat{T}_{\partial D}<\widehat{\tau}_{B(y,b_j|k|^{-\Theta})})\le 1-p,$$

then

$$\mathbf{P}_y(\widehat{T}_{D^c\cap B(y,(1/2)b_j|k|^{-\Theta})}<\widehat{\tau}_{B(y,b_j|k|^{-\Theta})})\le\mathbf{P}_y(\widehat{T}_{\partial D}<\widehat{\tau}_{B(y,b_j|k|^{-\Theta})})\le 1-p.$$

Since $b_j|k|^{-\Theta}\le r_1$, by Lemma 2.8 with $K:=\overline{B(0,1/2)}$ , we have

$$c_{15}^{-1}(b_j|k|^{-\Theta})^{2-d}\left(\inf_{w\in K}G^0_{B(0,1)}(0,w)\right)\mathrm{Cap}^{B(y,b_j|k|^{-\Theta})}(D^c\cap\overline{B(y,b_j|k|^{-\Theta}/2)})$$
$$\le\mathbf{P}_y(\widehat{T}_{D^c\cap B(y,(1/2)b_j|k|^{-\Theta})}<\widehat{\tau}_{B(y,b_j|k|^{-\Theta})})\le 1-p.$$

Thus,

$$\mathrm{Cap}^{B(y,b_j|k|^{-\Theta})}(D^c\cap\overline{B(y,b_j|k|^{-\Theta}/2)})\le c_{16}(1-p)b_j^{d-2}|k|^{-(d-2)\Theta}. \tag{4.13}$$

Using the facts that $D^c\cap D_2^{j,k}(y)\subset D^c\cap B(y,\frac{1}{2}b_j|k|^{-\Theta})$ and

$$|A\cap\overline{B(z,r/2)}|^{(d-2)/d}\le c_{17}\,\mathrm{Cap}^{B(z,r)}(A\cap\overline{B(z,r/2)}),\qquad z\in\mathbf{R}^d,$$

we have, from (4.13), that

$$\begin{aligned}|D^c \cap D_2^{j,k}(y)| &\le |D^c \cap B(y, \tfrac{1}{2} b_j |k|^{-\Theta})| \\ &\le c_{18} (\mathrm{Cap}^{B(y,b_j|k|^{-\Theta})} (D^c \cap \overline{B(y, b_j|k|^{-\Theta}/2)}))^{d/(d-2)} \\ &\le c_{19}(1-p)^{d/(d-2)} b_j^d |k|^{-d\Theta}.\end{aligned}$$

If we choose $p_1 \in (p_0, 1)$ and let $c_{20} := c_{19}(1-p_1)^{d/(d-2)} \max_{1 \le i \le m_0} b_i^d$ be such that

$$|D^c \cap D_2^{j,k}(y)| \le c_{20} |k|^{-d\Theta} = \tfrac{1}{2} |D_2^{j,k}(y)|,$$

then

$$|D \cap D_2^{j,k}(y)| > c_{20} |k|^{-d\Theta}/2.$$

Note that since $D$ is bounded, $D$ is an $L^d$-domain (a domain which can be locally represented as the region above the graph of an $L^d$ function). Now, we can follow the proof of Lemma 2.6 (with $p = d$ and $\Theta = r$ there) on the second half of page 417 in [2] (after rescaling) to get

$$(a_0 - y_d)/|k|^{-\Theta} \le c_{21} |k|^{\Theta(d-1)/d} |k|^{\Theta} = c_{21} |k|^{1-1/(4d)}. \tag{4.14}$$

Since $p_1 \in (p_0, 1)$, by Lemma 4.6,

$$\mathbf{P}_y(\widehat{T}_{M^{j,k}(y)} < \widehat{\tau}_{D_1^{j,k}(y)}) \ge \exp\left( -c_{22} \frac{8(a_0 - y_d)|k|^{\Theta}}{b_j} \right). \tag{4.15}$$

Using our (4.14)–(4.15) instead of (2.10)–(2.11) in [2], we can follow the argument in the proof of Lemma 2.6 after (2.11) in [2] (after rescaling) to conclude that $y \notin S_k$ if $-k$ is sufficiently large. Thus, we have proven the claim by contradiction. Moreover,

$$\mathbf{P}_y(\widehat{\tau}_{S_k} < \widehat{\tau}_{B(y, b_j|k|^{-\Theta})}) \ge \mathbf{P}_y(\widehat{\tau}_D < \widehat{\tau}_{B(y, b_j|k|^{-\Theta})}) > 1 - p_1, \qquad y \in S_k.$$

Thus, by Lemma 2.10, we have

$$\mathbf{E}_y[\widehat{\tau}_{S_k}] \le c_{23} \left( \max_{1 \le i \le m_0} b_i \right)^{-\Theta} |k|^{-2\Theta} = c_{24} |k|^{-1-\beta},$$

where $\beta = 1/(4d-2) > 0$. □

**5. Parabolic boundary Harnack principle and intrinsic ultracontractivity.** Throughout this section, we will assume that $D$ is one of the following types of bounded domains:

(a) a twisted Hölder domain of order $\alpha \in (1/3, 1]$;
(b) a uniformly Hölder domain of order $\alpha \in (0, 2)$;

(c) a bounded domain which can be locally represented as the region above the graph of a function.

Recall that $t_1$ is the constant from (3.1) and $\widehat{\tau}_B = \inf\{t > 0 : \widehat{Y}_t \notin B\}$. For any $\delta > 0$, we put $D_\delta := \{x \in D : \rho_D(x) < \delta\}$.

LEMMA 5.1. *There exist constants $c$, $R_1 > 0$ and a point $x_1$ in $D$ such that $B_1 := \overline{B(x_1, \frac{1}{2}R_1)} \subset D \setminus \overline{D_{(1/4)R_1}}$ and for every $R \le R_1$, $r^D(t,x,y) \wedge \widehat{r}^D(t,x,y) \ge cR^{-d}$ for all $x, y \in \overline{B(x_1, \frac{1}{2}R)}$ and $\frac{1}{3}t_1R^2 \le t \le t_1R^2$.*

PROOF. Choose $R_1 = R_1(D) \le \sqrt{t_0}$ and $x_1 \in D$ such that $B(x_1, R_1) \subset D$. We then apply Lemma 3.1 with $\delta = \frac{1}{3}$ and use the monotonicity of the density to get the desired assertion. $\square$

We fix $x_1$, $R_1$ and $B_1$ in the lemma above for the remainder of this section. Let $h_1(x) := G_D(x, x_1)$ and $h_2(x) := \widehat{G}_D(x, x_1)$. $h_1$ and $h_2$ are regular harmonic for $Y$ and $\widehat{Y}$ in $D \setminus B_1$, respectively. Moreover, by (2.4) and (2.6), $h_1$ and $h_2$ are bounded by $2^{k_0+1}$ for some $k_0 = k_0(R_1)$ on $D \setminus B_1$. Let $(\mathbf{P}^h_x, Y^D_t)$ and $(\mathbf{P}^h_x, \widehat{Y}^D_t)$ be the $h$-transforms of $(\mathbf{P}_x, Y^D_t)$ and $(\mathbf{P}_x, \widehat{Y}^D_t)$, respectively.

LEMMA 5.2. *For every $s > 0$, there exists a positive constant $\delta_0 = \delta_0(s) \le \frac{1}{4}R_1$ such that*

$$\left(\inf_{x\in D} \mathbf{P}^{h_1}_x\left(T_{D\setminus D_\delta} < \frac{s}{4}\right)\right) \wedge \left(\inf_{x\in D} \mathbf{P}^{h_2}_x\left(\widehat{T}_{D\setminus D_\delta} < \frac{s}{4}\right)\right) \ge \frac{1}{2}.$$

PROOF. For $k \le k_0$, let

$$V^\delta_k := \{x \in D_\delta : h_2(x) \le 2^{k+1}\}, \qquad U_k := \{x \in D \setminus B_1 : h_2(x) \le 2^{k+1}\}.$$

Clearly, $V^\delta_k \subset U_k$ for $\delta \le \frac{1}{4}R_1$. For each $k$, by (2.4) and (2.6), we have

$$\sup_{x\in D} \mathbf{E}_x[\widehat{\tau}_{V^\delta_k}] \le c \sup_{x\in D} \int_{V^\delta_k} \frac{dy}{|x-y|^{d-2}}$$

for some $c > 0$. So, $\sup_{x\in D} \mathbf{E}_x[\widehat{\tau}_{V^\delta_k}]$ goes to zero as $\delta \to 0$ by the uniform integrability of $|x-y|^{-d+2}$ over $D$. Note that $D \setminus B_1$ is also one of the types of domains we assumed at the beginning of this section. So, by Lemma 4.7,

$$\sum_{k=-\infty}^{k_0} \sup_{x\in D} \mathbf{E}_x[\widehat{\tau}_{U_k}] < \infty.$$

Thus, by the dominated convergence theorem, we have

$$\lim_{\delta\downarrow 0} \sum_{k=-\infty}^{k_0} \sup_{x\in D} \mathbf{E}_x[\widehat{\tau}_{V^\delta_k}] = 0. \tag{5.1}$$

On the other hand, since 1 is excessive for $\widehat{Y}^D$, it is easy to see that $(1/h_2(\widehat{Y}^D), \mathcal{F}_t)$ is a supermartingale with respect to $\mathbf{P}_x^{h_2}$, where $\mathcal{F}_t$ is the natural filtration of $\{\widehat{Y}^D\}$ (see, e.g., page 83 in [14]). Thus, with the same proof, one can see that the first inequality in equation (8) on page 179 of [8] is true. Thus, there exists $c_1$ independent of $h_2$ and $\delta$ such that

$$\sup_{x\in D} \mathbf{E}_x^{h_2}[\widehat{\tau}_{D_\delta}] \le c_1 \sum_{k=-\infty}^{k_0} \sup_{x\in D} \mathbf{E}_x[\widehat{\tau}_{V_k^\delta}]. \tag{5.2}$$

Combining (5.1)–(5.2), we have that for each $s>0$, there exists $\delta>0$ such that $\sup_{x\in D} \mathbf{E}_x^{h_2}[\widehat{\tau}_{D_\delta}] < s/8$. We can now apply Chebyshev's inequality to get

$$\mathbf{P}_x^{h_2}\left(\widehat{\tau}_{D_\delta} < \frac{s}{4}\right) \ge \frac{1}{2}.$$

On the other hand, using (2.2), (2.4) and (2.6), it is elementary to show that the strictly positive function $\widehat{G}_D(x,y)$ is $\infty$ if and only if $x=y\in D$, and for every $x\in D$, $\widehat{G}_D(x,\cdot)$ and $\widehat{G}_D(\cdot,x)$ are extended continuous in $D$ (see the proof in Theorem 2.6 in [16]). Thus, the condition (H) in [22] holds. Also, the strict positivity of $\widehat{G}_D(x,y)$ and Proposition 2.6 imply that the set $W$ on page 5 of [22] and the set $Z$ defined in [9] [equation (12) on page 179] are empty. Thus, by Theorem 2 in [22], for every $x\neq x_1$, the lifetime $\widehat{\zeta}^{h_2}$ of $\widehat{Y}^D$ is finite $\mathbf{P}_x^{h_2}$-a.s. and

$$\lim_{t\uparrow\zeta^{h_2}} \widehat{Y}_t^D = x_1, \qquad \mathbf{P}_x^{h_2}\text{-a.s.} \tag{5.3}$$

Thus, for $x\in D_\delta$, the conditioned process $\widehat{Y}^D$ with respect to $\mathbf{P}_x^{h_2}$ cannot be killed before hitting $D\setminus D_\delta$, due to the continuity of $\widehat{Y}^D$. Therefore, we have

$$\mathbf{P}_x^{h_2}\left(\widehat{T}_{D\setminus D_\delta} < \frac{s}{4}\right) = \mathbf{P}_x^{h_2}\left(\widehat{\tau}_{D_\delta} < \frac{s}{4}\right) \ge \frac{1}{2}. \qquad \square$$

For a parabolic function $g(t,x)$ in $\Omega=(T_1,T_2]\times D$ for $Y$ (resp. $\widehat{Y}$), let $(\mathbf{P}_{t,x}^g, Z_s^\Omega)$ [resp. $(\mathbf{P}_{t,x}^g, \widehat{Z}_s^\Omega)$] be the killed space-time process $(\mathbf{P}_{t,x}, Z_s^\Omega)$ [resp. $(\mathbf{P}_{t,x}, \widehat{Z}_s^\Omega)$] conditioned by $g$. For each $u>0$, we let

$$W_k = W_k(u) := \{(s,y)\in[u/2,u]\times D : 2^k \le g(s,x) \le 2^{k+1}\}$$

and

$$W^n = W^n(u) := \bigcup_{k=-\infty}^{n} W_k.$$

LEMMA 5.3. *For every $M>0$ and $u>0$, there exists $k_1=k_1(M,u,h_1,h_2,B_1)<-3$ such that for every positive parabolic function $g(t,x)$ in $(u/2,u]\times D$ for $Y$ (resp. $\widehat{Y}$),*

$$g(s,x)\ge Mh_1(x) \qquad [\text{resp. } g(s,x)\ge Mh_2(x)], \qquad (s,x)\in[u/2,u]\times(D\setminus B_1)$$

*implies*

$$\mathbf{E}^g_{u,x}[\tau_{W^{k_1}}]\le\frac{u}{8} \qquad \left(\text{resp. } \mathbf{E}^g_{u,x}[\widehat{\tau}_{W^{k_1}}]\le\frac{u}{8}\right), \qquad x\in D,$$

*where $\tau_{W^{k_1}}=\inf\{t>0: Z_t\notin W^{k_1}\}$ and $\widehat{\tau}_{W^{k_1}}=\inf\{t>0:\widehat{Z}_t\notin W^{k_1}\}$.*

PROOF. Let $m_1$ be the smallest integer greater than $\log_2 M$ and $U_k:=\{x\in D\setminus B_1: h_2(x)\le 2^{k+1}\}$ so that $W_k\subset U_{k+m_1}\times[u/2,u]$ for small $k$. By Lemma 4.7, we get, for small $n$,

$$\sum_{k=-\infty}^{n}\sup_{(s,y)\in W_k}\mathbf{E}_{s,y}[\widehat{\tau}_{W_k}]\le\sum_{k=-\infty}^{n}\sup_{(s,y)\in U_{k+m_1}}\mathbf{E}_{s,y}[\widehat{\tau}_{U_{k+m_1}}]<\infty. \tag{5.4}$$

Similarly to the argument in the proof of the previous lemma, using the estimates in [8], there exists $c_1$, independent of $g$, $n$ and $u$ and such that

$$\sup_{y\in D}\mathbf{E}^g_{u,y}[\widehat{\tau}_{W^n}]\le c_1\sum_{k=-\infty}^{n}\sup_{(s,y)\in W_k}\mathbf{E}_{s,y}[\widehat{\tau}_{W_k}]. \tag{5.5}$$

Combining (5.4)–(5.5), we have that, for small $n$,

$$\sup_{x\in D}\mathbf{E}^g_{u,x}[\widehat{\tau}_{W^n}]<\infty.$$

Now, choose $k_1=k_1(u)<0$ small so that

$$\sup_{x\in D}\mathbf{E}^g_{u,x}[\widehat{\tau}_{W^{k_1}}]<\frac{u}{8}.$$ □

The idea of the proof of the next lemma comes from the proof of Lemma 5.1 in [2]. We spell out the details for the reader's convenience.

LEMMA 5.4. *For every $u\in(0,\frac{1}{2}t_1R_1^2)$, there exists $c>0$ such that for all $x\in D$,*

$$\mathbf{P}_x(Y_u\in B_1,\tau_D>u)\ge c\mathbf{P}_x(\tau_D>u)$$

*and*

$$\mathbf{P}_x(\widehat{Y}_u\in B_1,\widehat{\tau}_D>u)\ge c\mathbf{P}_x(\widehat{\tau}_D>u).$$

PROOF. In this proof, for $A \subset [0,\infty)\times V$, $\widehat{T}_A$ will denote the first hitting time of $A$ for $\widehat{Z}_s$.

We fix $u \le \frac{1}{2}t_1R_1^2$, let $\delta_0 = \delta_0(u) \le \frac{1}{4}R_1$ be the constant from Lemma 5.2 and let $D_2 := D_{\delta_0}$. Note that $B_1 \subset D \setminus D_2$. Let $f_\varepsilon(x) = \varepsilon$ on $D\setminus B_1$ and 1 on $B_1$. Define a parabolic function $g_\varepsilon$ on $(0,\infty)\times D$ by

$$g_\varepsilon(t,x) := \int_D \widehat{r}^D(t,x,y) f_\varepsilon(y)\,dy = \mathbf{E}_x[f_\varepsilon(\widehat{Y}_t^D) : \widehat{Y}_t^D \in D], \qquad 0<\varepsilon<1.$$

Clearly,

$$\varepsilon \mathbf{P}_x(\widehat{\tau}_D > t) \le g_\varepsilon(t,x) \le \mathbf{P}_x(\widehat{\tau}_D > t). \tag{5.6}$$

We claim that there exists $c_1 > 0$ independent of $\varepsilon$ such that

$$g_\varepsilon(t,x) \ge c_1 h_2(x), \qquad (x,t) \in (D\setminus B_1)\times[u/2,u]. \tag{5.7}$$

First, we note that since $2u \le t_1R_1^2$, by Theorem 3.8 and a chain argument, we get

$$\begin{aligned} \inf_{(t,x)\in[u/4,u]\times(D\setminus D_2)} g_\varepsilon(t,x) &\ge c_1 g_\varepsilon(u/8,x_1) \\ &\ge c_1 \int_{B_1} \widehat{r}^D(u/8,x_1,y)\,dy = c_2 \end{aligned} \tag{5.8}$$

for some $c_2>0$. Let $h(t,x) := h_2(x)$ for $(t,x)\in[u/4,u]\times(D\setminus B_1)$. Since $h(t,x) \le 2^{k_0+1}$, by (5.8), we have

$$g_\varepsilon(t,x) \ge c_2 2^{-k_0-1} h(t,x), \qquad (t,y) \in [u/4,u]\times(D\setminus(D_2\cup B_1)).$$

Let $\Omega := (0,\infty)\times D$. For $(s,x)\in[u/2,u]\times D_2$,

$$\begin{aligned} g_\varepsilon(s,x) &\ge \mathbf{E}_{s,x}[g_\varepsilon(\widehat{Z}^\Omega_{\widehat{T}_{(0,\infty)\times(D\setminus D_2)}}) : \widehat{T}_{(0,\infty)\times(D\setminus D_2)} \le u/4] \\ &\ge c_2 2^{-k_0-1} h(s,x) \mathbf{P}^h_{s,x}(\widehat{T}_{(0,\infty)\times(D\setminus D_2)} \le u/4) \\ &= c_2 2^{-k_0-1} h_2(x) \mathbf{P}^h_{s,x}(\widehat{T}_{(0,\infty)\times(D\setminus D_2)} \le u/4) \\ &= c_2 2^{-k_0-1} h_2(x) \mathbf{P}^{h_2}_x(\widehat{T}_{D\setminus D_2} \le u/4), \end{aligned}$$

which is greater than or equal to $c_2 2^{-k_0-2} h_2(x)$, by Lemma 5.2. The claim is proved.

We now apply Lemma 5.3 to $g_\varepsilon(s,x)$ and get

$$\mathbf{E}^{g_\varepsilon}_{u,x}[\widehat{\tau}_{W^{k_1}}(\varepsilon)] \le \frac{u}{8}, \qquad x \in D. \tag{5.9}$$

Let $\varepsilon_1 := 2^{k_1-1} < \frac{1}{4}$, $g(s,x) := g_{\varepsilon_1}(s,x)$ and

$$E := W^{k_1} = \{(s,x)\in[u/2,u]\times D : g(s,x) \le 4\varepsilon_1\}.$$

By Chebyshev's inequality, from (5.9), we get

$$\mathbf{P}^g_{u,x}\left(\widehat{\tau}_E \le \frac{u}{4}\right) \ge \frac{1}{2}, \qquad x \in D. \tag{5.10}$$

Let $S_1$ be the first hitting time of $\partial(D \times [0,\infty))$ of $\widehat{Z}$. The conditioned process $(\mathbf{P}^g_{t,x}, \widehat{Z}^\Omega)$ cannot be killed before time $t$. In fact,

$$\begin{aligned}\mathbf{P}^g_{t,x}(\widehat{Z}^\Omega_{S_1-} \in \{0\} \times D) &= \mathbf{E}_{t,x}\left[\frac{g(\widehat{Z}^\Omega_{S_1-})}{g(t,x)} : Z^\Omega_{S_1-} \in \{0\} \times D\right] \\ &= \mathbf{E}_x\left[\frac{g(0,\widehat{Y}^D_t)}{g(t,x)} : \widehat{Y}^D_t \in D\right] \\ &= \frac{1}{g(t,x)}\mathbf{E}_x[f_{\varepsilon_1}(\widehat{Y}^D_t) : \widehat{Y}^D_t \in D] = 1.\end{aligned}$$

Thus, we get

$$\mathbf{P}^{g_\varepsilon}_{u,x}\left(\widehat{T}_{\partial_1 E} \le \frac{u}{4}\right) = \mathbf{P}^{g_\varepsilon}_{u,x}\left(\widehat{\tau}_E \le \frac{u}{4}\right) \ge \frac{1}{2}, \qquad x \in D, \tag{5.11}$$

where $\partial_1 E := \partial E \cap ((0,\infty) \times D)$.

Note that, by (5.6),

$$\begin{aligned}\mathbf{P}_x(\widehat{Y}_u \in B_1, \widehat{\tau}_D > u)/\mathbf{P}_x(\widehat{\tau}_D > u) &\ge \varepsilon_1 \mathbf{P}_x(\widehat{Y}_u \in B_1, \widehat{\tau}_D > u)/g(u,x) \\ &\ge \varepsilon_1 \mathbf{P}^g_{u,x}(\widehat{Z}^\Omega_{S_1-} \in \{0\} \times B_1).\end{aligned} \tag{5.12}$$

Thus, it is enough to bound $\mathbf{P}^g_{u,x}(\widehat{Z}^\Omega_{S_1-} \in \{0\} \times B_1)$. By the strong Markov property and (5.11),

$$\begin{aligned}\mathbf{P}^g_{u,x}(\widehat{Z}^\Omega_{S_1-} \in \{0\} \times B_1) &\ge \mathbf{P}^g_{u,x}\left(\widehat{Z}^\Omega_{S_1-} \in \{0\} \times B_1, \widehat{T}_{\partial_1 E} \le \frac{u}{4}\right) \\ &= \mathbf{E}^g_{u,x}\left[\mathbf{P}^g_{\widehat{Z}^\Omega_{\widehat{T}_{\partial_1 E}}}(\widehat{Z}^\Omega_{S_1-} \in \{0\} \times B_1) : \widehat{T}_{\partial_1 E} \le \frac{u}{4}\right] \\ &\ge \frac{1}{2} \inf_{(s,x)\in\partial_1 E} \mathbf{P}^g_{s,x}(\widehat{Z}^\Omega_{S_1-} \in \{0\} \times B_1).\end{aligned} \tag{5.13}$$

Since $g = 4\varepsilon_1$ on $\partial_1 E$ by the continuity of $g$, for $(s,x) \in \partial_1 E$,

$$\begin{aligned}4\varepsilon_1 &= \int_D \widehat{r}^D(s,x,y) f_{\varepsilon_1}(y)\,dy \\ &= \mathbf{P}_x(\widehat{Y}^D_s \in B_1) + \varepsilon_1 \mathbf{P}_x(\widehat{Y}^D_s \in D \setminus B_1) \\ &= \mathbf{P}_{s,x}(\widehat{Z}_{S_1} \in \{0\} \times B_1) + \varepsilon_1 \mathbf{P}_{s,x}(\widehat{Z}_{S_1} \in \{0\} \times (D \setminus B_1)) \\ &\le \mathbf{P}_{s,x}(\widehat{Z}_{S_1} \in \{0\} \times B_1) + \varepsilon_1.\end{aligned}$$

Thus,

$$\mathbf{P}_{s,x}(\widehat{Z}_{S_1} \in \{0\} \times B_1) \geq 3\varepsilon_1.$$

Since $\mathbf{P}^g_{s,x}(\widehat{Z}^{\Omega}_{S_1-} \in \{0\} \times D) = 1$, applying the above inequality, we get

$$\begin{aligned} &\mathbf{P}^g_{s,x}(\widehat{Z}^{\Omega}_{S_1-} \in \{0\} \times B_1) \\ &\quad = \frac{1}{4\varepsilon_1}\mathbf{E}_{s,x}[g(\widehat{Z}^{\Omega}_{S_1-}); \widehat{Z}^{\Omega}_{S_1-} \in \{0\} \times B_1] \\ &\quad = \frac{1}{4\varepsilon_1}\mathbf{P}_{s,x}(\widehat{Z}_{S_1} \in \{0\} \times B_1) \geq \frac{3}{4} > 0, \qquad (s,x) \in \partial_1 E. \end{aligned} \tag{5.14}$$

Combining (5.11)–(5.14), the proof is completed. □

Let $p(t,x,y) := r^D(t,x,y)/H(y)$. Recall that $H(y) = \int_V G(x,y)\,dx$ and $\xi(dy) = H(y)\,dy$. For any $t > 0$, define

$$P^D_t f(x) := \int_D r^D(t,x,y) f(y)\,dy = \int_D p(t,x,y) f(y) \xi(dy)$$

and

$$\widehat{P}^D_t f(x) := \int_D \widehat{r}^D(t,x,y) f(y)\,dy = \int_D p(t,y,x) f(y) \xi(dy).$$

By definition, we have

$$\int_D f(x) P^D_t g(x) \xi(dx) = \int_D g(x) \widehat{P}^D_t f(x) \xi(dx).$$

It is easy to check that $\{P_t\}$ and $\{\widehat{P}_t\}$ are both strongly continuous contraction semigroups in $L^2(D, \xi(dx))$. We will use $\mathcal{L}$ and $\widehat{\mathcal{L}}$ to denote the $L^2(D,\xi(dx))$-infinitesimal generators of $\{P^D_t\}$ and $\{\widehat{P}^D_t\}$, respectively.

LEMMA 5.5.

(1)

$$\frac{p(t,x,y)}{p(t,x,z)} \geq c_1 \frac{p(t,w,y)}{p(t,w,z)} \qquad \forall w,x,y,z \in D$$

*implies that for every $s > t$ and $w,x,y,z \in D$,*

$$\frac{p(s,y,x)}{p(s,z,x)} \geq c_1 \frac{p(t,y,w)}{p(t,z,w)} \quad \textit{and} \quad \frac{p(s,x,y)}{p(s,x,z)} \leq c_1^{-1} \frac{p(t,w,y)}{p(t,w,z)}.$$

(2)

$$\frac{p(t,y,x)}{p(t,z,x)} \geq c_2 \frac{p(t,y,v)}{p(t,z,v)} \qquad \forall v,x,y,z \in D$$

*implies that for every $s > t$ and $v, x, y, z \in D$,*

$$\frac{p(s,x,y)}{p(s,x,z)} \geq c_2 \frac{p(t,v,y)}{p(t,v,z)} \quad \textit{and} \quad \frac{p(s,y,x)}{p(s,z,x)} \leq c_2^{-1} \frac{p(t,y,v)}{p(t,z,v)}.$$

PROOF. We give the proof of (2) only. The proof of (1) is similar. Since

$$p(t,w,y) \geq c_2 \frac{p(t,w,z)}{p(t,v,z)} p(t,v,y) \qquad \forall w,x,y,z \in D,$$

we get

$$\begin{aligned} p(s,x,y) &= \int_D p(s-t,x,w) p(t,w,y) \xi(dw) \\ &\geq c_2 \frac{p(t,v,y)}{p(t,v,z)} \int_D p(s-t,x,w) p(t,w,z) \xi(dw) \\ &= \frac{p(t,v,y)}{p(t,v,z)} p(s,x,z). \end{aligned}$$

On the other hand, since

$$p(t,y,w) \leq c_2^{-1} \frac{p(t,y,v)}{p(t,z,v)} p(t,z,w) \qquad \forall w,x,y,z \in D,$$

we get

$$\begin{aligned} p(s,y,x) &= \int_D p(t,y,w) p(s-t,w,x) \xi(dw) \\ &\leq c_2^{-1} \frac{p(t,y,v)}{p(t,z,v)} \int_D p(t,z,w) p(s-t,w,x) \xi(dw) \\ &= \frac{p(t,y,v)}{p(t,z,v)} p(s,z,x). \end{aligned}$$

□

THEOREM 5.6. *For each $u \in (0, \frac{1}{2} t_1 R_1^2)$, there exists $c = c(D,u) > 0$ such that*

$$\frac{p(t,x,y)}{p(t,x,z)} \geq c \frac{p(s,v,y)}{p(s,v,z)}, \qquad \frac{p(t,y,x)}{p(t,z,x)} \geq c \frac{p(s,y,v)}{p(s,z,v)} \tag{5.15}$$

*for every $s, t \geq u$ and $v, x, y, z \in D$.*

PROOF. Let $\tau_1 := \inf\{t > 0 : Y_t \notin D\}$, $\tau_2 := \inf\{t > 0 : \widehat{Y}_t \notin D\}$, $\varphi_1(x) := \mathbf{P}_x(\tau_1 > u/3)$ and $\varphi_2(y) := \mathbf{P}_y(\tau_2 > u/3)$. By (2.1) with $T = \frac{1}{2} t_1 R_1^2$, there exists $c_1 > 0$ such that

$$p(u,x,y) = \int_D p\left(\frac{u}{3}, x, z\right) \int_D p\left(\frac{u}{3}, z, w\right) p\left(\frac{u}{3}, w, y\right) \xi(dw) \xi(dz)$$

$$\begin{aligned}&\le c_1 u^{-d/2}\int_D p\left(\frac{u}{3},x,z\right)\xi(dz)\int_D p\left(\frac{u}{3},w,y\right)\xi(dw)\\&= c_1 u^{-d/2}\varphi_1(x)\varphi_2(y).\end{aligned}$$

For the lower bound, we use Lemmas 5.1 and 5.4, and get

$$\begin{aligned}p(u,x,y)&\ge \int_{B_1} p\left(\frac{u}{3},x,z\right)\int_{B_1} p\left(\frac{u}{3},z,w\right)p\left(\frac{u}{3},w,y\right)\xi(dw)\xi(dz)\\&\ge c_2 u^{-d/2}\int_{B_1} p\left(\frac{u}{3},x,z\right)\xi(dz)\int_{B_1} p\left(\frac{u}{3},w,y\right)\xi(dw)\\&= c_2 u^{-d/2}\mathbf{P}_x(Y_{u/3}\in B_1,\tau_1>u)\mathbf{P}_y(\widehat{Y}_{u/3}\in B_1,\tau_2>u)\\&\ge c_3 u^{-d/2}\varphi_1(x)\varphi_2(y)\end{aligned}$$

for some positive constants $c_2$ and $c_3$. Thus, both inequalities in (5.15) are true for $s=t=u\le \frac{1}{2}t_1R_1^2$. We now apply Lemma 5.5(1)–(2) and get, for $s>u$ and $v,x,y,z\in D$,

$$\frac{p(s,y,x)}{p(s,z,x)}\ge c_4\frac{p(u,y,v)}{p(u,z,v)},\qquad \frac{p(s,x,y)}{p(s,x,z)}\le c_4^{-1}\frac{p(u,v,y)}{p(u,v,z)} \tag{5.16}$$

and

$$\frac{p(s,x,y)}{p(s,x,z)}\ge c_4\frac{p(u,v,y)}{p(u,v,z)},\qquad \frac{p(s,y,x)}{p(s,z,x)}\le c_4^{-1}\frac{p(u,y,v)}{p(u,z,v)}. \tag{5.17}$$

Thus, both inequalities in (5.15) are true for $s>t=u$. Moreover, combining (5.16)–(5.17), both inequalities in (5.15) are also true for $t=s>u$. Again applying Lemma 5.5(1)–(2), we complete the proof. □

By (2.6), we have proven the parabolic boundary Harnack principle for $Y^D$.

COROLLARY 5.7. *For each positive $u\in(0,\frac{1}{2}t_1R_1^2)$, there exists $c=c(D,u)>0$ such that*

$$\frac{r^D(t,x,y)}{r^D(t,x,z)}\ge c\frac{r^D(s,w,y)}{r^D(s,w,z)},\qquad \frac{r^D(t,y,x)}{r^D(t,z,x)}\ge c\frac{r^D(s,y,w)}{r^D(s,z,w)}$$

*for every $s,t\ge u$ and $w,x,y,z\in D$.*

Since, for each $t>0$, $p(t,x,y)$ is bounded in $D\times D$, it follows from Jentzsch's theorem (Theorem V.6.6 on page 337 of [26]) that the common value $\lambda_0:=\sup\mathrm{Re}(\sigma(\mathcal{L}))=\sup\mathrm{Re}(\sigma(\widehat{\mathcal{L}}))$ is an eigenvalue of multiplicity 1 for both $\mathcal{L}$ and $\widehat{\mathcal{L}}$, that an eigenfunction $\phi_0$ of $\mathcal{L}$ associated with $\lambda_0$ can be chosen to be strictly positive with $\|\phi_0\|_{L^2(D,\xi(dx))}=1$ and an eigenfunction $\psi_0$ of $\widehat{\mathcal{L}}$ associated with $\lambda_0$ can be chosen to be strictly positive with $\|\psi_0\|_{L^2(D,\xi(dx))}=1$.

DEFINITION 5.8. The semigroups $\{P_t^D\}$ and $\{\widehat{P}_t^D\}$ are said to be intrinsically ultracontractive if for any $t>0$, there exists a constant $c_t>0$ such that

$$p(t,x,y) \le c_t \phi_0(x)\psi_0(y) \qquad \forall (x,y)\in D\times D.$$

Now, the next theorem, which is the main result of this paper, can be easily proven from Lemma 5.4 and the continuity of $\phi_0$ and $\psi_0$. But we give the proof that Theorem 5.6 implies the intrinsic ultracontractivity.

THEOREM 5.9. *The semigroups $\{P_t^D\}$ and $\{\widehat{P}_t^D\}$ are intrinsically ultracontractive. Moreover, for any $t>0$, there exists a constant $c_t>0$ such that*

$$c_t^{-1}\phi_0(x)\psi_0(y) \le p(t,x,y) \le c_t\phi_0(x)\psi_0(y) \qquad \forall (x,y)\in D\times D. \tag{5.18}$$

PROOF. Integrating both sides of (5.15) with respect to $y$ over $D$ for $t=s=u\le \frac{1}{2}t_1R_1^2$, we get

$$\frac{p(t,x,z)}{\int_D p(t,x,y)\xi(dy)} \le c_t \frac{p(t,w,z)}{\int_D p(t,w,y)\xi(dy)} \tag{5.19}$$

and

$$\frac{p(t,z,x)}{\int_D p(t,y,x)\xi(dy)} \le c_t \frac{p(t,z,w)}{\int_D p(t,y,w)\xi(dy)} \tag{5.20}$$

for all $w,x,z\in D$. We fix $x_0\in D$. The above (5.20) implies that for any positive function $f$ and $z\in D$,

$$\begin{aligned}
P_t^D f(z) &= \int_D p(t,z,x)f(x)\xi(dx)\\
&\le c_t\Big(\int_D p(t,y,x_0)\xi(dy)\Big)^{-1}\int_D\int_D p(t,y,x)\xi(dy)p(t,z,x_0)f(x)\xi(dx)\\
&= c_t\frac{p(t,z,x_0)}{\int_D p(t,y,x_0)\xi(dy)}\int_D\int_D p(t,y,x)\xi(dy)f(x)\xi(dx)\\
&= c_t\frac{p(t,z,x_0)}{\int_D p(t,y,x_0)\xi(dy)}\int_D P_t^D f(y)\xi(dy).
\end{aligned}$$

Similarly, (5.20) also implies the lower bound

$$P_t^D f(z) \ge c_t^{-1}\frac{p(t,z,x_0)}{\int_D p(t,y,x_0)\xi(dy)}\int_D P_t^D f(y)\xi(dy), \qquad z\in D.$$

Using (5.19), we also get the corresponding result for $\widehat{P}_t^D$. Thus, we have, for all $z, w \in D$,

$$c_t^{-1}\frac{p(t,z,x_0)}{\int_D p(t,y,x_0)\xi(dy)} \le \frac{P_t^D f(z)}{\int_D P_t^D f(y)\xi(dy)} \le c_t\frac{p(t,z,x_0)}{\int_D p(t,y,x_0)\xi(dy)} \tag{5.21}$$

and

$$c_t^{-1}\frac{p(t,x_0,w)}{\int_D p(t,x_0,y)\xi(dy)} \le \frac{\widehat{P}_t^D f(w)}{\int_D \widehat{P}_t^D f(y)\xi(dy)} \le c_t\frac{p(t,x_0,w)}{\int_D p(t,x_0,y)\xi(dy)}. \tag{5.22}$$

Applying (5.21) to $\phi_0$ and a sequence of functions approaching the point mass at $w$ appropriately, we get that for any $z, w \in D$,

$$c_t^{-2}\phi_0(z) \le \frac{p(t,z,w)}{\int_D p(t,y,w)\xi(dy)} \le c_t^2\phi_0(z),$$

which implies that

$$c_t^{-4}\frac{\phi_0(z)}{\phi_0(x_0)} \le \frac{p(t,z,w)}{p(t,x_0,w)} \le c_t^4\frac{\phi_0(z)}{\phi_0(x_0)}, \qquad z,w \in D. \tag{5.23}$$

Similarly, applying (5.22) to $\psi_0$ and a sequence of functions approaching the point mass at $z$, we get that for any $z, w \in D$,

$$c_t^{-4}\frac{\psi_0(w)}{\psi_0(x_0)} \le \frac{p(t,z,w)}{p(t,z,x_0)} \le c_t^4\frac{\psi_0(w)}{\psi_0(x_0)}. \tag{5.24}$$

Thus, combining (5.23)–(5.24), we conclude that for any $t \le \frac{1}{2}t_1R_1^2$ and any $z, w \in D$,

$$\begin{aligned} p(t,z,w) &= p(t,x_0,x_0)\frac{p(t,x_0,w)}{p(t,x_0,x_0)}\frac{p(t,z,w)}{p(t,x_0,w)} \\ &\le c_t^8 p(t,x_0,x_0)\frac{\phi_0(z)\psi_0(w)}{\phi_0(x_0)\psi_0(x_0)}. \end{aligned}$$

Let $T := \frac{1}{2}t_1R_1^2$. Since

$$\begin{aligned} p(s,x,y) &= \int_D p(T,x,z)p(s-T,z,y)\xi(dz) \\ &\le c_T^8 c_{s-T}^8 \phi_0(x)\psi_0(y)\int_D \phi_0(z)\psi_0(z)\xi(dz) \\ &\le c_T^8 c_{s-T}^8 \phi_0(x)\psi_0(y), \qquad s \in (T,2T], \end{aligned}$$

we can easily get the intrinsic ultracontractivity by induction. The fact that intrinsic ultracontractivity implies the lower bound is proved in [18] (Proposition 2.5 in [18]). □

Let

$$\begin{aligned} \phi(x) &:= \phi_0(x) \Big/ \int_D \phi_0(y)^2 \, dy, \\ \psi(x) &:= \psi_0(x) H(x) \Big/ \int_D \psi_0(y)^2 H(y)^2 \, dy. \end{aligned} \tag{5.25}$$

Note that $0 < \int_D \psi_0(y)^2 H(y)^2 \, dy < \infty$ because of (2.6). Since

$$e^{\lambda_0 t} \phi_0(x) = \int_D p(t,x,y) \phi_0(y) \xi(dy) = \int_D r^D(t,x,y) \phi_0(y) \, dy$$

and

$$\begin{aligned} e^{\lambda_0 t} \psi_0(x) H(x) &= H(x) \int_D p(t,y,x) \psi_0(y) \xi(dy) \\ &= \int_D r^D(t,y,x) \psi_0(y) H(y) \, dy, \end{aligned}$$

we have

$$\begin{aligned} e^{\lambda_0 t} \phi(x) &= \int_D r^D(t,x,y) \phi(y) \, dy, \\ e^{\lambda_0 t} \psi(x) &= \int_D r^D(t,y,x) \psi(y) \, dy. \end{aligned} \tag{5.26}$$

We say that the common value $e^{\lambda_0 t}$ is an eigenvalue for $r^D(t,x,y)$ and the pair $(\phi, \psi)$ are the corresponding eigenfunctions if (5.26) is true and if $\phi$ and $\psi$ are strictly positive with $\|\phi\|_{L^2(D,dx)} = 1$ and $\|\psi\|_{L^2(D,dx)} = 1$. So, the intrinsic ultracontractivity of $\{P_t^D\}$ and $\{\widehat{P}_t^D\}$ can be rephrased as follows.

Corollary 5.10. *For any $t > 0$, there exists a constant $c_t > 0$ such that*

$$c_t^{-1} \phi(x) \psi(y) \le r^D(t,x,y) \le c_t \phi(x) \psi(y) \qquad \forall (x,y) \in D \times D. \tag{5.27}$$

Proof. This is clear from (5.18) and (5.25). □

Applying Theorem 2.7 of [18], we have the following.

Theorem 5.11. *There exist positive constants $c$ and $a$ such that for every $(t,x,y) \in (1,\infty) \times D \times D$,*

$$\left| \left( e^{-\lambda_0 t} \int_D \phi_0(z) \psi_0(z) \xi(dz) \right) \frac{r^D(t,x,y)}{\phi_0(x) \psi_0(y) H(y)} - 1 \right| \le c e^{-at}. \tag{5.28}$$

We are going to use $SH^+$ to denote families of nonnegative superharmonic functions of $Y$ in $D$. For any $h \in SH^+$, we use $\mathbf{P}_x^h$ to denote the law of the $h$-conditioned diffusion process $Y^D$ and $\mathbf{E}_x^h$ to denote the expectation with respect to $\mathbf{P}_x^h$. Let $\zeta^h$ be the lifetime of the $h$-conditioned diffusion process $Y^D$.

In [18], the bound for the lifetime of the conditioned $Y^D$ is proved using Theorem 5.11.

Theorem 5.12 (Theorem 2.8 in [18]).
(1)

$$\sup_{x \in D, h \in SH^+} \mathbf{E}_x^h[\zeta^h] < \infty.$$

(2) *For any* $h \in SH^+$, *we have*

$$\lim_{t \uparrow \infty} e^{-\lambda_0 t} \mathbf{P}_x^h(\zeta^h > t) = \frac{\phi_0(x)}{h(x)} \int_D \psi_0(y) h(y) \xi(dy) \Big/ \int_D \phi_0(y) \psi_0(y) \xi(dy).$$

*In particular,*

$$\lim_{t \uparrow \infty} \frac{1}{t} \log \mathbf{P}_x^h(\zeta^h > t) = \lambda_0.$$

**Acknowledgments.** The authors are grateful to the referees for their valuable comments.

Department of Mathematics
Seoul National University
Seoul 151-742
Republic of Korea
E-mail: pkim@snu.ac.kr

Department of Mathematics
University of Illinois
Urbana, Illinois 61801
USA
E-mail: rsong@math.uiuc.edu